\documentclass[12pt,twoside]{article}

\usepackage{amsmath,amssymb,amsthm,amscd,graphicx,color,accents}

\usepackage[colorlinks=true,linkcolor=blue]{hyperref}
\usepackage{caption,lipsum}

\numberwithin{equation}{section}

\theoremstyle{plain}
\newtheorem{theorem}{Theorem}[section]
\newtheorem{proposition}[theorem]{Proposition}

\newtheorem{corollary}[theorem]{Corollary}
\newtheorem{remark}[theorem]{Remark}
\newcommand{\Pro}{\par\noindent\rm\textsc{Proof. }}

\newcommand{\sbt}{\, \begin{picture}(-1,1)(-1,-2)\circle*{2}\end{picture}\ }
\newcommand{\sbts}{\, \begin{picture}(-1,1)(-1,-1.5)\circle*{2}\end{picture}\ }

\def\e{{{\hskip 1pt}\rm{e}{\hskip 1pt}}}

\def\bk{{\boldsymbol{k}}}
\def\bJ{{\boldsymbol{J}}}
\def\bK{{\boldsymbol{K}}}
\def\bN{{\boldsymbol{N}}}

\def\bzero{\boldsymbol{0}}
\def\dd{{{\hskip1pt}\rm{d}}}
\def\im{{\rm{i}}\hskip1pt}

\def\P{{\mathcal P}}
\def\Q{{\mathcal Q}}
\def\V{{\mathcal V}}
\def\E{{\rm{E}}\hskip 1pt}
\def\C{\mathbb{C}}
\def\N{\mathbb{N}}
\def\R{\mathbb{R}}

\def\re{{\rm{Re}}}

\def\SX{{\Sigma_{X}}}
\def\SY{{\Sigma_{Y}}}
\def\SXY{{\Sigma_{XY}}}
\def\SYX{{\Sigma_{Y\!X}}}
\def\LXY{{\Lambda_{XY}}}

\begin{document}

\title{\bf 
Distance Correlation Coefficients for Lancaster Distributions
}

\author{
{Johannes Dueck,}\thanks{Institute of Applied Mathematics, University of Heidelberg, Im Neuenheimer Feld 294, 69120 Heidelberg, Germany.} 
\ \ {Dominic Edelmann,}\thanks{Division of Biostatistics, German Cancer Research Center (DKFZ), Im Neuenheimer Feld 280, 69120 Heidelberg, Germany.} 
\ {\ and Donald Richards}\thanks{Department of Statistics, Pennsylvania State University, University Park, PA 16802, U.S.A.
\endgraf
\ $^\ddag$Corresponding author; e-mail address: richards@stat.psu.edu}}

\bigskip


\maketitle

\begin{abstract} 
We consider the problem of calculating distance correlation coefficients between random vectors whose joint distributions belong to the class of Lancaster distributions.  We derive under mild convergence conditions a general series representation for the distance covariance for these distributions.  To illustrate the general theory, we apply the series representation to derive explicit expressions for the distance covariance and distance correlation coefficients for the bivariate normal distribution and its generalizations of Lancaster type, the multivariate normal distributions, and the bivariate gamma, Poisson, and negative binomial distributions which are of Lancaster type.  



\smallskip
\noindent
{\em Running head}: Distance Correlation and Lancaster Distributions.
\end{abstract}

\medskip

\section{Introduction}
\label{sec:introduction}

The concepts of distance covariance and distance correlation, introduced by Sz\'ekely, et al. \cite{szekely2009brownian,szekely2007measuring}, have been shown to be widely applicable for measuring dependence between collections of random variables.  As examples of the ubiquity of distance correlation methods, we note the results on distance correlation given recently by: Sz\'ekely, et al. \cite{rizzo2010disco,szekely2012uniqueness,szekely2013distance,szekely2014partial,szekely2007measuring}, on statistical inference; Sejdinovic, et al. \cite{sejdinovic2013equivalence}, on machine learning; Kong, et al. \cite{kong2012using}, on familial relationships and mortality; Zhou \cite{zhou2012measuring}, on nonlinear time series; Lyons \cite{lyons2013distance}, on abstract metric spaces; Mart\'inez-Gomez, et al.~\cite{martinez2014distance} and Richards, et al.~\cite{richards2014interpreting}, on large astrophysical databases; Dueck, et al. \cite{dueck2014affinely}, on high-dimensional inference and the analysis of wind data; and Dueck, et al. \cite{dueck2015generalization}, on a connection with singular integrals on Euclidean spaces.  

A result which is of fundamental importance in distance correlation theory is the explicit formula for the empirical distance correlation coefficient \cite[pp.~2773-2774]{szekely2007measuring}.  By combining that explicit formula with the fast algorithm of Huo and Sz\'ekely \cite{huo2015fast}, it becomes straightforward to apply distance correlation methods to real-world data sets.  

On the other hand, the calculation of population distance correlation coefficients remains an intractable problem generally.  Sz\'ekely, et al. \cite[pp.~2785-2786]{szekely2007measuring} calculated the distance correlation coefficient for the bivariate normal distribution; Dueck, et al. \cite[Appendix]{dueck2012affinely} extended that result to the general multivariate normal distribution; and Dueck, et al. \cite{dueck2014affinely} calculated the affinely invariant distance correlation coefficient for the multivariate normal distribution.  Otherwise, no such results are yet available for any other distribution.  Hence, the state of distance correlation theory hitherto is that the empirical coefficients can be calculated readily but the opposite holds for their population counterparts, generally.  Consequently, it was not possible to calculate distance correlation coefficients explicitly for given nonnormal distributions in terms of the usual parameters that parametrize these distributions, or to ascertain for nonnormal distributions any analogs of the limit theorems derived by Dueck, et al. \cite[Section 4]{dueck2014affinely}.  

We describe in detail the difficulties arising in attempts to calculate the population distance correlation coefficients.  Let $p$ and $q$ be positive integers.  For column vectors $s \in \R^p$ and $t \in \R^q$, denote by $\|s\|$ and $\|t\|$ the standard Euclidean norms on the corresponding spaces; thus, if $s = (s_1,\ldots,s_p)^\top$ then $\|s\| = (s_1^2+\cdots+s_p^2)^{1/2}$, and similarly for $\|t\|$.  Given vectors $u$ and $v$ of the same dimension, we let $\langle u,v \rangle$ be the standard Euclidean scalar product of $u$ and $v$.  For jointly distributed random vectors $(X,Y) \in \R^p \times \R^q$ and non-random vectors $(s,t) \in \R^p \times \R^q$, let 
$$
\psi_{X,Y}(s,t) = \E \exp \! \big[\,\im \langle s,X \rangle + \im \langle t,Y \rangle \big],
$$
$\im = \sqrt{-1}$, be the joint characteristic function of $(X,Y)$, and let $\psi_X(s) = \psi_{X,Y}(s,0)$ and $\psi_Y(t) = \psi_{X,Y}(0,t)$ be the corresponding marginal characteristic functions.  For any $z \in \C$, let $|z|^2$ denote the squared modulus of $z$; also, we use the notation 
\begin{equation}
\label{eq:cp}
\gamma_p = \frac{\pi^{(p+1)/2}}{\Gamma\big((p+1)/2\big)}. 
\end{equation}
In the case of distributions with finite first moments, Sz\'ekely, et al. \cite[p.~2772]{szekely2007measuring} defined $\V(X,Y)$, the {\em distance covariance} between $X$ and $Y$, to be the positive square-root of 
\begin{equation}
\label{eq:dcov}
\V^2(X,Y) = \frac{1}{\gamma_p \gamma_q} \int_{\R^{p+q}} 
\frac{|\psi_{X,Y}(s,t) - \psi_X(s)\psi_Y(t)|^2}
{\|s\|^{p+1} \, \|t\|^{q+1}} \, \dd s \, \dd t
\end{equation}
and they defined the {\em distance correlation coefficient} between $X$ and $Y$ as 
\begin{equation}
\label{eq:dcor}
\mathcal{R}(X,Y) = \frac{\V(X,Y)}{\sqrt{\V(X,X)\V(Y,Y)}}
\end{equation}
if both $\V(X,X)$ and $\V(Y,Y)$ are strictly positive, and otherwise to be zero \cite[p.~2773]{szekely2007measuring}.  For distributions with finite first moments we have $0 \le \mathcal{R}(X,Y) \le 1$, and $\mathcal{R}(X,Y) = 0$ if and only if $X$ and $Y$ are mutually independent.  

For given random vectors $X$ and $Y$, the fundamental obstacle in calculating the population distance correlation coefficient (\ref{eq:dcor}) is the computation of the singular integral (\ref{eq:dcov}).  In particular, the singular nature of the integrand precludes evaluation of the integral by expanding the numerator, $|\psi_{X,Y}(s,t) - \psi_X(s)\psi_Y(t)|^2$, and subsequent term-by-term integration of each of the resulting three terms.  

In this paper, we calculate the distance correlation coefficients for pairs $(X,Y)$ of random vectors whose joint distributions are in the class of {\em Lancaster distributions}, a class of probability distributions made prominent by Lancaster \cite{lancaster1958structure,lancaster1969} and Sarmanov \cite{sarmanov1966generalized}.  The distribution functions of the Lancaster family are well-known to have attractive expansions in terms of certain orthogonal functions (Koudou \cite{koudou1998lancaster}; Diaconis, et al.~\cite{diaconis2008gibbs}).  By applying those expansions, we obtain explicit expressions for the distance covariance and distance correlation coefficients.  

Consequently, we derive under mild convergence conditions a general formula for the distance covariance for the Lancaster distributions.  We apply the general formula to obtain explicit expressions for the distance covariance and distance correlation for the bivariate normal distributions and some of its generalizations, for the multivariate normal distributions, and for bivariate gamma, Poisson, and negative binomial distributions.  We remark that explicit results can also be obtained for other Lancaster-type expansions obtained by Bar-Lev, et al. \cite{bar1994diagnonal}; however, we will omit the details for other cases because the formulas derived here are entirely representative of other cases.

\section{The Lancaster distributions}
\label{sec:SL}

To recapitulate the class of Lancaster distributions we generally follow the standard notation in that area, as given by Koudou \cite{koudou1996probabilites,koudou1998lancaster}; cf., Lancaster \cite{lancaster1969}, Pommeret \cite{pommeret2004characterization}, or Diaconis, et al. \cite[Section 6]{diaconis2008gibbs}.  

Let $(\mathcal{X},\mu)$ and $(\mathcal{Y},\nu)$ be locally compact, separable probability spaces, such that $L^2(\mu)$ and $L^2(\nu)$ are separable.  Let $\sigma$, a probability measure on $\mathcal{X}\times\mathcal{Y}$, have marginal distributions $\mu$ and $\nu$; then there exist functions $K_\sigma$ and $L_\sigma$ such that 
$$
\sigma(\dd x,\dd y) = K_\sigma(x,\dd y) \mu(\dd x) = L_\sigma(\dd x,y) \nu(\dd y).
$$
We note that $K_\sigma$ and $L_\sigma$ represent the conditional distributions of $Y$ given $X=x$, and $X$ given $Y=y$, respectively.  

Let $\mathcal{C}$ denote a countable index set with a zero element, denoted by $0$.  Let $\{P_n: n \in \mathcal{C}\}$ and $\{Q_n: n \in \mathcal{C}\}$ be sequences of functions on $\mathcal{X}$ and $\mathcal{Y}$ which form orthonormal bases for the separable Hilbert spaces $L^2(\mu)$ and $L^2(\nu)$, respectively.  We assume, by convention, that $P_0 \equiv 1$ and $Q_0 \equiv 1$.  

Since the tensor product Hilbert space $L^2(\mu \otimes \nu) \equiv L^2(\mu) \otimes L^2(\nu)$ is separable there holds, for $\sigma \in L^2(\mu \otimes \nu)$, the expansion 
\begin{equation}
\label{sigma_expansion}
\sigma(\dd x,\dd y)  = \sum_{m \in \mathcal{C}} \sum_{n \in \mathcal{C}} \rho_{m,n} P_m(x) Q_n(y) \, \mu(\dd x) \, \nu(\dd y),
\end{equation}
$(x,y) \in \mathcal{X} \times \mathcal{Y}$.  Letting $\delta_{m,n}$ denote Kronecker's delta, the probability measure $\sigma$ is called a {\textit{Lancaster distribution}} if there exists a nonnegative sequence $\{\rho_n: n \in \mathcal{C}\}$ such that
$$
\int P_m(x) \, Q_n(y) \, \sigma(\dd x,\dd y) = \rho_m \, \delta_{m,n}
$$
for all $m,n \in \mathcal{C}$; in particular, $\rho_0 = 1$.  The sequence $\{\rho_n: n \in \mathcal{C}\}$ is called a {\textit{Lancaster sequence}}, and the expansion (\ref{sigma_expansion}) reduces to 
$$
\sigma(\dd x,\dd y)  = \sum_{n \in \mathcal{C}} \rho_n P_n(x) Q_n(y) \mu(\dd x) \nu(\dd y).
$$
Koudou \cite[pp.~255--256]{koudou1996probabilites} characterized the Lancaster sequences $\{\rho_n: n \in \mathcal{C}\}$ such that the associated probability distribution $\sigma$ is absolutely continuous with respect to $\mu \otimes \nu$ and has Radon-Nikodym derivative 
$$
\frac{\sigma(\dd x,\dd y)}{\mu(\dd x) \, \nu(\dd y)} = \sum_{n \in \mathcal{C}} \rho_n \, P_n(x) \, Q_n(y) \ \in \ L^2(\mu \otimes \nu),
$$
$(x,y) \in \mathcal{X} \times \mathcal{Y}$.  

In the sequel, we consider the case in which $\mathcal{X} = \mathbb{R}^p$ and $\mathcal{Y} = \mathbb{R}^q$ and the underlying random vectors $X \in \mathbb{R}^p$ and $Y \in \mathbb{R}^q$ have joint distribution $\sigma$ and marginal distributions $\mu$ and $\nu$, respectively.  We assume that $\mu$, $\nu$, and $\sigma$ are absolutely continuous with respect to Lebesgue measure or counting measure on the respective sample spaces and we denote their corresponding probability density functions by $\phi_X$, $\phi_Y$, and $\phi_{X,Y}$, respectively.  This yields the expansion, 
\begin{equation} 
\label{lancaster.expansion}
\phi_{X,Y}(x,y) = \phi_X(x) \, \phi_Y(y) \sum_{n \in \mathcal{C}} \rho_n \, P_n(x) \, Q_n(y).
\end{equation}
We will refer to (\ref{lancaster.expansion}) as the {\textit{Lancaster expansion}} of the joint density function $\phi_{X,Y}$.

\section{Examples of Lancaster expansions}
\label{Lancaster.Examples}

In this section, we provide examples of Lancaster expansions (\ref{lancaster.expansion}) for the bivariate normal distribution and some of its generalizations, the multivariate normal distributions, and the bivariate gamma, Poisson, and negative binomial distributions.  In the sequel, we denote by $\N_0$ the set of nonnegative integers.

\subsection{The bivariate normal distribution and some of its generalizations}
\label{lancaster_bvn}

Let $(X,Y)$ follow a bivariate normal distribution with mean vector $0$ and covariance matrix 
$$
\Sigma = \begin{pmatrix} 1 & \rho \\ \rho & 1 \end{pmatrix},
$$
denoted by $(X,Y) \sim \mathcal{N}_2(0,\Sigma)$.  The joint probability density function of $(X,Y)$ is 
$$
\phi_{X,Y}(x,y;\rho) = \frac{1}{2\pi}(1-\rho^2)^{-\tfrac12}\,\exp\left(-\frac{x^2+y^2-2\rho\,x\,y}{2(1-\rho^2)}\right),
$$
$x, y \in \R$, and the marginal density functions are given by 
$$
\phi_X(x) = \phi_Y(x) = \frac{1}{\sqrt{2\pi}}\exp\left(-\tfrac12 x^2\right).
$$ 
In this case, the index set $\mathcal{C}$ is $\N_0$.  For $n \in \N_0$, let 
$$
H_n(x) = (-1)^n \exp\big(\tfrac12 x^2) \Big(\frac{\dd}{\dd x}\Big)^n \exp\left(-\tfrac12 x^2\right),
$$
$x \in \R$, denote the $n$th Hermite polynomial, $n=0,1,2,\ldots$.  It is well-known that the polynomials $\{H_n: n \in \N_0\}$ are orthogonal with respect to the standard normal distribution and form a complete orthogonal basis for the Hilbert space $L^2(X)$.  Also, the Lancaster expansion of $\phi_{X,Y}$ is given by the classical formula of Mehler, which states that, for $x, y \in \R$, 
\begin{align}
\label{SL-biv-normal}
\phi_{X,Y}(x,y;\rho) = \phi_X(x) \, \phi_Y(y) \sum_{n=0}^\infty \frac{\rho^n}{n!} \, H_n(x) \, H_n(y),
\end{align}  
and this series converges absolutely for all $x, y \in \R$.  

\medskip

We remark that there are numerous extensions of Mehler's formula which represent Lancaster-type expansions for generalizations of the bivariate normal distribution.  Sarmanov and Bratoeva \cite{sarmanov1967probabilistic} consider series expansions of the form 
\begin{equation}
\label{sarmanov_bratoeva}
\phi_{X,Y}(x,y) = \phi_X(x) \, \phi_Y(y) \sum_{n=0}^\infty \frac{\rho_n}{n!} \, H_n(x) H_n(y),
\end{equation}
$x, y \in \R$, where the sequence of real numbers $\{\rho_n:n=0,1,2,\ldots\}$ satisfies $\sum_{n=0}^\infty \rho_n^2 < \infty$.  Sarmanov and Bratoeva proved that for the expansion (\ref{sarmanov_bratoeva}) to be nonnegative, and therefore to be a valid probability density function, it is necessary and sufficient that the sequence $\{\rho_n\}$ be the moment sequence of a random variable $\xi$ supported on the interval $[-1,1]$.  

An example of this generalization is the case in which $\phi_{X,Y}$ is a mixture of bivariate normal densities; put in the formula 
\begin{align}
\phi_{X,Y}(x,y) &= \tfrac12\big[\phi_{X,Y}(x,y;\rho)+\phi_{X,Y}(x,y;-\rho)\big] \nonumber \\
&= \phi_X(x) \phi_Y(y)\sum_{n=0}^\infty \frac{\rho^{2n}}{(2n)!} \, H_{2n}(x) H_{2n}(y) ;\label{eq:bvn_uncor}
\end{align}
this corresponds to the case in (\ref{sarmanov_bratoeva}) in which 
$$
\rho_n = \begin{cases} \rho^n, & n \hbox{ even} \\ 0, & n \hbox{ odd} \end{cases}
$$
This mixture density also provides an example of a distribution for which the Pearson correlation coefficent is zero whereas the distance correlation is positive.  


\subsection{The multivariate normal distribution}
\label{lancaster_mvn}

Let $X \in \R^p$ and $Y \in \R^q$ be random vectors such that $(X,Y) \sim \mathcal{N}_{p+q}(0,\Sigma)$, a $(p+q)$-dimensional multivariate normal distribution with mean vector $0$ and positive definite covariance matrix 
\begin{equation}
\label{Sigmamatrix}
\Sigma = \left( \begin{matrix} \SX & \SXY \\ \SYX & \SY \end{matrix} \right)
\end{equation}
where $\SX$, $\SY$, and $\SXY = \SYX'$ are $p \times p$, $q \times q$ and $p \times q$ matrices, respectively.  We denote by $\phi_{X,Y}$ the joint probability density function of $(X,Y)$, and by $\phi_X$ and $\phi_Y$ the marginal density functions of $X$ and $Y$, respectively.  

We now describe the Lancaster expansion of $\phi_{X,Y}$, a result derived in \cite{withers2010expansions}.  In this case, the index set $\mathcal{C}$ is $\N_0^{p \times q}$, the set of $p \times q$ matrices with nonnegative integer entries.  

For a matrix of summation indices $\bN = (N_{rc}) \in \N_0^{p \times q}$, define $\bN! = \prod_{r=1}^p \prod_{c=1}^q N_{rc}!$.  For $r=1,\ldots,p$, let 
$$
\bN_{r \sbt} = \sum_{c=1}^q N_{rc}
$$
and set $\bN_{* \sbt} = (\bN_{1 \sbt},\ldots,\bN_{p \sbt})$.  Similarly, for each $c=1,\ldots,q$, define 
$$
\bN_{\sbt c} = \sum_{r=1}^p N_{rc}
$$
and set $\bN_{\sbt *} = (\bN_{\sbt 1},\ldots,\bN_{\sbt q})$.  Further, we define 
$$
\bN_{\sbt\sbt} = \sum_{r=1}^p \sum_{c=1}^q N_{rc},
$$
and note that $\bN_{\sbt\sbt} = \sum_{r=1}^p \bN_{r \sbt} = \sum_{c=1}^q \bN_{\sbt c}$.  

Denoting by $(\SXY)_{rc}$ the $(r,c)$th entry of $\SXY$, we also define 
$$
\SXY^\bN = \prod_{r=1}^p \prod_{c=1}^q [(\SXY)_{rc}]^{N_{rc}}.
$$

We now introduce the multivariate Hermite polynomials.  For any $p \in \N$, $\bk = (k_1,\ldots,k_p) \in \N_0^p$, and $x = (x_1,\ldots,x_p) \in \R^p$, define $x^\bk = x_1^{k_1} \cdots x_p^{k_p}$ and define the differential operator, 
$$
\left(-\frac{\partial}{\partial x}\right)^\bk = \left(-\frac{\partial}{\partial x_1}\right)^{k_1} \cdots \left(-\frac{\partial}{\partial x_p}\right)^{k_p}.
$$ 
The $\bk$th {\it multivariate Hermite polynomial} with respect to the marginal density function $\phi_X$ is defined as 
\begin{equation}
\label{mvhermitepoly}
H_{\boldsymbol{k}}(x;\SX) = \frac{1}{\phi_X(x)} \, \left(-\frac{\partial}{\partial x}\right)^\bk \phi_X(x).
\end{equation}
The Lancaster expansion of the multivariate normal density function $\phi_{X,Y}$ is given by the generalized Mehler formula \cite{withers2010expansions}: 
\begin{equation}
\label{lancaster_multinormal}
\phi_{X,Y}(x,y) = \phi_{X}(x) \, \phi_{Y}(y) \sum_{\bN \in \N_0^{p \times q}} \frac{\SXY^\bN}{\bN!} \, H_{\bN_{* \sbts}}(x;\SX) \, H_{\bN_{{\sbts }*}}(y;\SY),
\end{equation}
with absolute convergence for all $x \in \R^p$, $y \in \R^q$.  

To calculate the affinely invariant distance correlation coefficient between $X$ and $Y$, as defined by Dueck, et al. (2014), we need the Lancaster expansion of the joint density function of the standardized random vectors $\widetilde{X} = \SX^{-1/2} \, X$ and $\widetilde{Y} = \SY^{-1/2} \, Y$.  It is straightforward to verify that $(\widetilde{X},\widetilde{Y}) \sim \mathcal{N}_{p+q}(0,\Lambda)$ where 
\begin{equation}
\label{Lambdamatrix}
\Lambda = \left(\begin{matrix} I_p & \LXY \\ \LXY' & I_q \end{matrix}\right)
\end{equation}
with $\LXY = \SX^{-1/2} \, \SXY \, \SY^{-1/2}$, and then we deduce from (\ref{lancaster_multinormal}) that the Lancaster expansion for $(\widetilde{X},\widetilde{Y})$ is 
\begin{align} 
\label{multi.normal.lancaster}
\phi_{\widetilde{X},\widetilde{Y}}(x,y) = \phi_{\widetilde{X}}(x) \, \phi_{\widetilde{Y}}(y) \sum_{\bN \in \N_0^{p \times q}} \frac{\LXY^\bN}{\bN!} \, H_{\bN_{* \sbts}}(x;I_p) \, H_{\bN_{\sbts *}}(y;I_q).
\end{align}

\subsection{The bivariate gamma distribution}
\label{lancaster_bvgamma}

The Lancaster expansion for a bivariate gamma distribution, which was derived by Sarmanov \cite{sarmanov1970gamma,sarmanov1970approximate}, can be stated as follows (see Kotz, et al. \cite[pp. 437--438]{kotz2000continuous}).  

For $\alpha > -1$ and $n \in \N_0$, the classical {\it Laguerre polynomial} is defined by 
\begin{equation}
\label{laguerre_explicit}
\begin{aligned}
\widetilde{L}_n^{(\alpha)}(x) &= \frac{1}{n!} x^{-\alpha} \exp(x) \Big(\frac{\dd}{\dd x}\Big)^n x^{n+\alpha} \exp(-x) \\
&= \frac{(\alpha+1)_n}{n!} \sum_{j=0}^n \frac{(-n)_j}{(\alpha+1)_j} \, \frac{x^j}{j!},
\end{aligned}
\end{equation}
$x > 0$, where 
$$
(\alpha)_n = \frac{\Gamma(\alpha+n)}{\Gamma(\alpha)} = \alpha(\alpha+1)\cdots(\alpha+n-1),
$$
$n=0,1,2,\ldots$, denotes the {\it rising factorial}.  By standardizing the classical Laguerre polynomial, we obtain the orthonormal version \cite{griffiths1969canonical},
$$
L_n^{(\alpha)}(x) = \left( \frac{(\alpha+1)_n}{n!}\right)^{-1/2}  \,  \widetilde{L}_n^{(\alpha)}(x)  = \left(\frac{(\alpha+1)_n}{n!} \right)^{1/2}\, \sum_{j=0}^n \frac{(-n)_j}{(\alpha+1)_j} \, \frac{x^j}{j!}.
$$
	
Let $\lambda \in (0,1)$, and let $\alpha$ and $\beta$ satisfy $\alpha \ge \beta > 0$.  Sarmanov \cite{sarmanov1970gamma,sarmanov1970approximate} derived for certain bivariate gamma random variables $(X,Y)$ the joint probability density function, 
\begin{equation}
\label{density.gamma}
\phi_{X,Y}(x,y) = \phi_X(x) \, \phi_Y(y) \sum_{n=0}^\infty a_n L_n^{(\alpha-1)}(x) \, L_n^{(\beta-1)}(y),
\end{equation}
$x, y > 0$, where 
\begin{equation}
\label{akcoefficients}
a_n = \left[\frac{(\beta)_n}{(\alpha)_n}\right]^{1/2} \, \lambda^{n},
\end{equation}
$n=0,1,2,\ldots$.  The corresponding marginal density functions are 
$$
\phi_X(x) = \frac{1}{\Gamma(\alpha)} \, x^{\alpha-1} \, \exp(-x)
$$
and 
$$
\phi_Y(y) = \frac{1}{\Gamma(\beta)} \, y^{\beta-1} \, \exp(-y),
$$
which we recognize as the density functions of one-dimensional gamma random variables with index parameters $\alpha$ and $\beta$, respectively.  

We remark that if $\alpha = \beta$ then the density function (\ref{density.gamma}) reduces to the Kibble-Moran bivariate gamma density function, $\hbox{Corr}(X,Y) = \lambda$ \cite[pp. 436--437]{kotz2000continuous}, and (\ref{density.gamma}) represents the Lancaster expansion for $(X,Y)$.  On the other hand, if $\alpha \neq \beta$ then $\hbox{Corr}(X,Y) \not\equiv \lambda$.  

More generally, Griffiths \cite{griffiths1969canonical} showed that a series expansion of the form 
\begin{equation}
\label{density.gamma.gen}
\phi_{X,Y}(x,y) = \phi_X(x) \, \phi_Y(y) \sum_{n=0}^\infty \rho_n L_n^{(\alpha-1)}(x) \, L_n^{(\beta-1)}(y)
\end{equation}
represents a valid bivariate probability density if and only if 
\begin{equation}
\label{griff.coef}
\rho_n = \left[\frac{(\beta)_n}{(\alpha)_n}\right]^{1/2} \, \lambda_n,
\end{equation}	
where $\lambda_n$ is the moment sequence of a random variable $\xi$ concentrated on $[0,1]$.

\subsection{The bivariate Poisson distribution}
\label{lancaster_bvpoisson}

For $a > 0$ and $x, n \in \N_0$, let 
\begin{equation}
\label{poisson_charlier}
C_n(x;a) = \Big(\frac{a^n}{n!}\Big)^{1/2} \sum_{k=0}^n (-1)^k \binom{n}{k} \binom{x}{k} \, \frac{k!}{a^k}
\end{equation}
denote the Poisson-Charlier polynomial of degree $n$.  For $\lambda \in [0,1]$, Koudou \cite[Section 5]{koudou1998lancaster} (cf., Bar-Lev, et al. \cite{bar1994diagnonal}, Pommeret \cite{pommeret2004characterization}) showed that there exists a bivariate random vector $(X,Y)$ with probability density function 
\begin{equation}
\label{bivariatepoisson_pdf}
\phi_{X,Y}(x,y) = \phi_X(x) \, \phi_Y(y) \sum_{n=0}^\infty \lambda^n\, C_n(x;a) \, C_n(y;a),
\end{equation}
$x, y \in \N_0$.  The corresponding marginal density functions $\phi_X$ and $\phi_Y$ are given by 
$$
\phi_X(k) = \phi_Y(k) = \frac{a^k \, \exp(-a)}{k!},
$$
$k \in \N_0$, so that $X$ and $Y$ are distributed marginally according to a Poisson distribution with parameter $a$.  The series (\ref{bivariatepoisson_pdf}) is an expansion of Lancaster type, a special case of (\ref{lancaster.expansion}), and the resulting distribution is called a bivariate Poisson distribution.

\subsection{The bivariate negative binomial distribution}

The orthonormal polynomials for the classical univariate negative binomial distribution are the (normalized) Meixner polynomials, given by
\begin{equation}
\label{meixner}
M_n^{\beta,c}(x) = \Big(\frac{c^n\,(\beta)_n}{n!}\Big)^{1/2} \sum_{k=0}^n \frac{(-n)_k \, (-x)_k}{(\beta)_k\,k!} \left(1-c^{-1}\right)^k,
\end{equation}
for $\beta > 0$, $0 < c < 1$, and $x \in \mathbb{N}$.  Koudou \cite[Section 6]{koudou1998lancaster} showed, by an approach similar to that used for the bivariate Poisson distribution, that there exists a bivariate random variable $(X,Y)$, with identical marginal negative binomial densities, 
$$
\phi_X(x) = \phi_Y(x) = (1-c)^\beta\, \frac{c^x\,(\beta)_x}{x!},
$$
$x \in \N_0$, and with joint probability density function, 
\begin{equation}
\label{bivariatenb_pdf}
\phi_{X,Y}(x,y) = \phi_X(x)\,\phi_Y(y) \sum_{n=0}^\infty \lambda^n\, M_n^{\beta,c}(x)\,M_n^{\beta,c}(y),
\end{equation}
where $x, y \in \N_0$, and $0 \leq \lambda < 1$.  The expansion (\ref{bivariatenb_pdf}) represents a Lancaster expansion of the joint density function.

\section{Distance correlation coefficients for Lancaster distributions}
\label{sec:DCforL}

In this section, we derive a general series expression for the distance correlation coefficients for Lancaster distributions with density functions of the form (\ref{lancaster.expansion}).  
For a joint density function $\phi_{X,Y}$ given by (\ref{lancaster.expansion}) and $n \in \mathcal{C}$, we introduce the notation 
\begin{equation}
\label{lnx}
\P_n(s) = \E \exp(\im\,\langle s,X\rangle) \, P_n(X),
\end{equation}
$s \in \R^p$, and 
\begin{equation}
\label{lny}
\Q_n(t) = \E \exp(\im\,\langle t,Y\rangle) \, Q_n(Y),
\end{equation}
$t \in \R^q$.  To verify that the expectation in (\ref{lnx}) converges absolutely for all $s \in \R^p$, we apply the Cauchy-Schwarz inequality to obtain 
$$
\E |\exp(\im\langle s,X\rangle) P_n(X)| \le \left(\E |\exp(\im\langle s,X\rangle)|^2\right)^{1/2} \cdot \left(\E |P_n(X)|^2\right)^{1/2} = 1,
$$
because $\{P_n: n \in \mathcal{C}\}$ is an orthonormal basis for the Hilbert space $L^2(\mu)$.  In particular, 
\begin{align*}
|\P_n(s)| \le \E |\exp(\im\langle s,X\rangle) P_n(X)| \le 1,
\end{align*}
for all $s \in \R^p$ and, similarly, $|\Q_n(t)| \le 1$ for all $t \in \R^q$.  

In the following result, we will use the notation 
$$
\mathcal{A}_{j,k} = \int_{\R^p} \P_j(s) \, \P_k(-s) \, \frac{\dd s}{\|s\|^{p+1}}
$$
and 
$$
\mathcal{B}_{j,k} = \int_{\R^q} \Q_j(t) \, \Q_k(-t) \, \frac{\dd t}{\|t\|^{q+1}},
$$
$j, k \in \mathcal{C}$, whenever these integrals converge absolutely.  

We now state the main result.  

\begin{theorem}
\label{Dcov_Lancaster.theorem}
Suppose that the random vectors $X \in \mathbb{R}^p$ and $Y \in \mathbb{R}^q$ have the joint probability density function (\ref{lancaster.expansion}).  Then, 
\begin{equation}
\label{Dcov_Lancaster}
\V^2(X,Y) = \frac{1}{\gamma_p \gamma_q} \, \sum_{j \in \mathcal{C}, j \neq 0} \, \sum_{k \in \mathcal{C}, k \neq 0} \rho_j \, \rho_k \, \mathcal{A}_{j,k} \, \mathcal{B}_{j,k},
\end{equation}
whenever the sum converges absolutely.
\end{theorem}

\Pro 
Rewriting the Lancaster expansion (\ref{lancaster.expansion}) in the form, 
$$
\phi_{X,Y}(x,y) - \phi_X(x) \, \phi_Y(y) = \phi_X(x) \, \phi_Y(y) \sum_{n \in \mathcal{C}, n \neq 0} \rho_n P_n(x) \, Q_n(y),
$$
and taking Fourier transforms on both sides of this identity, we obtain for all $s \in \R^p$ and $t \in \R^q$ the expansion 
\begin{equation}
\label{A1} 
\psi_{X,Y}(s,t) - \psi_X(s) \, \psi_Y(t) = \sum_{n \in \mathcal{C}, n \neq 0} \rho_n  \P_n(s) \,  \Q_n(t).
\end{equation}
This identity is valid subject to the requirement that we may interchange summation and integration, which is justified by the assumption that the sum in the final result converges absolutely.  Using (\ref{A1}) we deduce that 
\begin{align*}
|\psi_{X,Y}(s,t) - \psi_X(s)\psi_Y(t)|^2 &= \big(\,\psi_{X,Y}(s,t) - \psi_X(s)\psi_Y(t)\,\big) \big(\,\overline{\psi_{X,Y}(s,t) - \psi_X(s)\psi_Y(t)}\,\big) \\
&= \sum_{j \in \mathcal{C}, j \neq 0} \, \sum_{k \in \mathcal{C}, k \neq 0} \rho_j \, \rho_k \, \P_j(s) \, \P_k(-s) \, \Q_j(t) \, \Q_k(-t).
\end{align*}
Next, we integrate this expansion with respect to the measures $\dd s/\|s\|^{p+1}$ and $\dd t/\|t\|^{q+1}$; this requires that we again interchange summation and integration which, by assumption, we are able to do.  On carrying through these procedures, we obtain (\ref{Dcov_Lancaster}).
$\qed$

\medskip

\begin{remark}
\label{Dcov_Lancaster_alpha_remark}
{\rm 
We note that Theorem \ref{Dcov_Lancaster.theorem} can be extended to the more general $\alpha$-distance covariance and correlation measures treated by Sz\'ekely, et al. \cite[p.~2784]{szekely2007measuring}.  For $\alpha \in (0,2)$, define 
$$
\gamma_{p,\alpha} = \frac{2 \, \pi^{p/2} \, \Gamma\big(1-\tfrac12\alpha\big)}{\alpha 2^\alpha \Gamma\big(\tfrac12(p+\alpha)\big)} 
$$
and let 
$$
\V^2_\alpha(X,Y) = \frac{1}{\gamma_{p,\alpha} \gamma_{q,\alpha}} \int_{\R^{p+q}} 
\frac{|\psi_{X,Y}(s,t) - \psi_X(s)\psi_Y(t)|^2}
{\|s\|^{p+\alpha} \, \|t\|^{q+\alpha}} \, \dd s \, \dd t
$$
be the $\alpha$-distance covariance between the random vectors $X$ and $Y$.  Further, define 
$$
\mathcal{A}_{j,k}(\alpha) = \int_{\R^p} \P_j(s) \, \P_k(-s) \, \frac{\dd s}{\|s\|^{p+\alpha}}
$$
and 
$$
\mathcal{B}_{j,k}(\alpha) = \int_{\R^q} \Q_j(t) \, \Q_k(-t) \, \frac{\dd t}{\|t\|^{q+\alpha}},
$$
$j, k \in \mathcal{C}$, whenever these integrals converge absolutely.  Then, the extension of Theorem \ref{Dcov_Lancaster.theorem} to the $\alpha$-distance covariance measures is that 
$$
\V^2_\alpha(X,Y) = \frac{1}{\gamma_{p,\alpha} \gamma_{q,\alpha}} \, \sum_{j \in \mathcal{C}, j \neq 0} \, \sum_{k \in \mathcal{C}, k \neq 0} \rho_j \, \rho_k \, \mathcal{A}_{j,k}(\alpha) \, \mathcal{B}_{j,k}(\alpha),
$$
whenever the sum converges absolutely.  The proof of this result is similar to the proof of Theorem \ref{Dcov_Lancaster.theorem}.  
}\end{remark}

\section{Examples}
\label{sec:examples}

In this section, we establish the versatility of Theorem \ref{Dcov_Lancaster.theorem} by applying it to compute the distance correlation coefficients for the bivariate normal, multivariate normal, and bivariate gamma, Poisson, and negative binomial distributions.  We verify for each example the absolute convergence of the series resulting from Theorem \ref{Dcov_Lancaster.theorem}, for that convergence property cannot be obtained in general from the general theorem.  In developing each example, we retain the corresponding notation in Section \ref{Lancaster.Examples}.  
We also remark that the crucial singular integral in the theory of distance correlation \cite{dueck2015generalization,szekely2007measuring} is evaluated in terms of the gamma function; however, even slight generalizations of that integral can be evaluated only in terms of the Gaussian or the confluent hypergeometric series; this explains the appearance of those series in the ensuing examples.

\subsection{The bivariate normal distribution and some of its generalizations}
\label{sub_dcovbvn}

In the sequel, we use the standard double-factorial notation, 
$$
n!! = n(n-2)(n-4)\cdots 
= \begin{cases} 
1, & \hbox{ if } n = -1, 0 \\
n(n-2)(n-4) \cdots 1, & \hbox{ if } n = 1, 3, 5, 7, \ldots \\
n(n-2)(n-4) \cdots 2, & \hbox{ if } n = 2, 4, 6, 8, \ldots 
\end{cases}
$$

\begin{proposition}
Let $(X,Y) \sim \mathcal{N}_2(0,\Sigma)$, a bivariate normal distribution with correlation coefficient $\rho$.  Then, 
\begin{equation}
\label{dcov_bvn}
\V^2(X,Y) = 4\pi^{-1} \sum_{\ell=1}^\infty \frac{\left((2\ell-3)!!\right)^2}{(2\ell)!} \big(1 - 2^{-(2\ell-1)}\big) \rho^{2\ell},
\end{equation}
and this series converges absolutely for all $\rho \in (-1,1)$.  
\end{proposition}    

\Pro 
Starting with the Lancaster expansion of the bivariate normal density function, as given in (\ref{SL-biv-normal}), and using the definitions of $\P_n$ and $\Q_n$ in (\ref{lnx}) and (\ref{lny}), respectively, we obtain by substitution and integration-by-parts, 
\begin{align*}
\P_n(s) = \Q_n(s) &= \int_{-\infty}^\infty \exp(\im sx) \frac{1}{\sqrt{2\pi}} \exp\left(-\tfrac{1}{2}x^2\right) \, H_n(x) \dd x \\
&= (\im s) ^n \exp\left(-\tfrac{1}{2}s^2\right),
\end{align*}
$s \in \R$.  Therefore, 
\begin{eqnarray*}
\mathcal{A}_{j,k} = \mathcal{B}_{j,k} 
&=& (-1)^{k} \, \im^{j+k} \int_{-\infty}^\infty  \,s^{j+k-2} \, \exp(-s^2) \, \dd s, \\
&=& \begin{cases}
(-1)^{k} \, \im^{j+k} \, \pi^{1/2} \left(\tfrac12\right)^{(j+k-2)/2}\,(j+k-3)!!, & \hbox{if } j+k \hbox{ is even} \\
0, & \hbox{otherwise }
\end{cases}
\end{eqnarray*}
since the latter integral is a moment of the $\mathcal{N}(0,\tfrac12)$ distribution.  By Theorem \ref{Dcov_Lancaster.theorem}, we obtain 
$$
\V^2(X,Y) = \frac{4}{\pi} \sum_{ \substack{j, \, k \, > \, 0 \\ j+k \ {\rm{even}}} }\frac{\rho^{j+k}}{j!\,k!}\,\left(\tfrac12\right)^{j+k}\,\big((j+k-3)!!\big)^2.
$$
Setting $j+k=2\ell$ with $\ell \ge 1$, the double series reduces to 
\begin{eqnarray*}
\V^2(X,Y) &=& \frac{4}{\pi} \sum_{\ell=1}^\infty \rho^{2\ell} (\tfrac12)^{2\ell} ((2\ell-3)!!)^2 \sum_{\substack{j, k \ge 1 \\ j+k = 2\ell}} \frac{1}{j!\,k!} \\
&=& \frac{4}{\pi} \sum_{\ell=1}^\infty \rho^{2\ell} (\tfrac12)^{2\ell} \frac{((2\ell-3)!!)^2}{(2\ell)!} \sum_{j=1}^{2\ell-1} \frac{(2\ell)!}{j!\,(2\ell-j)!} \\
&=& \frac{4}{\pi} \sum_{\ell=1}^\infty \rho^{2\ell} (\tfrac12)^{2\ell} \frac{((2\ell-3)!!)^2}{(2\ell)!} (2^{2\ell} - 2),
\end{eqnarray*}
which is the same as (\ref{dcov_bvn}).  

The absolute convergence of (\ref{dcov_bvn}) can be verified by comparison with a geometric series.  Moreover, it can be shown that the series reduces to the explicit formula, 
\begin{equation} 
\label{eq:dcov_bvn_ex}
\frac{4}{\pi} \, \Big( \rho \sin^{-1} \rho + \sqrt{1-\rho^2} - \rho \sin^{-1}  \rho/2 - \sqrt{4-\rho^2} +1 \Big)
\end{equation}  
which is identical with the result obtained by Sz\'ekely, et al. \cite[pp.~2785-2786]{szekely2007measuring}. 
$\qed$

\medskip

Having obtained $\mathcal{V}(X,Y)$, we let $\rho \to 1-$ to obtain the distance variances $\mathcal{V}(X,X)$ and $\mathcal{V}(Y,Y)$; here, we are applying a well-known result that if $(X,Y) \sim \mathcal{N}_2(0,\Sigma)$ where $\hbox{Var}(X) = \hbox{Var}(Y)$ and $\rho = 1$ then $X = Y$, almost surely.  By applying properties of Gauss' hypergeometric series, as was done by Dueck, et al. \cite[p.~2318]{dueck2014affinely}, we obtain 
$$
\mathcal{V}^2(X,X) = \mathcal{V}^2(Y,Y) = \frac{4}{3} - \frac{4(\sqrt{3}-1)}{\pi}.
$$

It is straightforward to extend the above results to generalizations of the type given in equation (\ref{sarmanov_bratoeva}).

\begin{corollary}
\label{cor:dcov_bvn_gen}
Let $(X,Y)$ be a bivariate random variable distributed according to a density function as given in (\ref{sarmanov_bratoeva}). Then
\begin{equation}
\label{dcov_bvn_gen}
\V^2(X,Y) = \frac{4}{\pi} \sum_{ \substack{j, \, k \, > \, 0 \\ j+k \ {\rm{even}}} }\frac{\rho_j \, \rho_k}{j!\,k!}\,\left(\tfrac12\right)^{j+k}\,\big((j+k-3)!!\big)^2.
\end{equation}
\end{corollary}

For the example given in (\ref{eq:bvn_uncor}), the series expansion in (\ref{dcov_bvn_gen}) reduces to an explicit formula similar to (\ref{eq:dcov_bvn_ex}).

\begin{corollary}
\label{cor:dcov_bvn_uncor}
Let $(X,Y)$ be a bivariate random variable distributed according to a density function as given in (\ref{sarmanov_bratoeva}) with $\rho_n = \rho^n$ for $n$ even and $\rho_n = 0$ for $n$ odd. Then 
\begin{equation}
\label{dcov_bvn_uncor}
\V^2(X,Y) = \frac{4}{\pi} \, \Big( \frac{\rho}{2} \sin^{-1} \rho + \tfrac12 \, \sqrt{1-\rho^2} - \rho \sin^{-1}  \rho/2 - \sqrt{4-\rho^2} +\frac{3}{2} \Big).
\end{equation}
\end{corollary}

\Pro  
Proceeding as in the case of the bivariate normal distribution, we obtain 
\begin{eqnarray*}
\V^2(X,Y) &=& \frac{4}{\pi} \sum_{\ell=1}^\infty \rho^{2\ell} (\tfrac12)^{2\ell} ((2\ell-3)!!)^2 \sum_{\substack{j, k \ge 1 \\ j+k = \ell}} \frac{1}{(2j)!\,(2k)!} \\
&=& \frac{4}{\pi} \sum_{\ell=1}^\infty \rho^{2\ell} (\tfrac12)^{2\ell} \frac{((2\ell-3)!!)^2}{(2\ell)!} \sum_{j=1}^{\ell-1} \frac{(2\ell)!}{(2j)!\,(2\ell-2j)!} \\
&=& \frac{4}{\pi} \sum_{\ell=1}^\infty \rho^{2\ell} (\tfrac12)^{2\ell} \frac{((2\ell-3)!!)^2}{(2\ell)!} (2^{2\ell-1} - 2). \\
\end{eqnarray*}
Using the standard notation, ${}_2F_1$, for Gauss' hypergeometric function we see that 
\begin{eqnarray*}
\V^2(X,Y) &=&\frac{4}{\pi} \sum_{\ell=1}^\infty \rho^{2\ell} (\tfrac12)^{2\ell} \frac{(2^\ell \,(-\tfrac12)_\ell)^2}{2^{2\ell} \, \ell! \,  (\tfrac12)_\ell} (2^{2\ell-1} - 2) \\
&=& \frac{4}{\pi} \Big[ \tfrac12 \, \big({}_2F_1(-\tfrac12,-\tfrac12;\tfrac12;\rho^2) - 1\big) - 2 \big( {}_2F_1(-\tfrac12,-\tfrac12;\tfrac12;\tfrac14 \rho^2) - 1 \big) \Big].
\end{eqnarray*}
It is well-known (see Andrews, Askey, and Roy~\cite[, pages 64 and 94]{andrews1999special}) that 
$$
{}_2F_1(-\tfrac12,-\tfrac12;\tfrac12;\rho^2) = \rho \sin^{-1}\rho + (1-\rho^2)^{1/2}.
$$
On applying this formula to the above expression, we obtain (\ref{dcov_bvn_uncor}).  
\qed

\subsection{The multivariate normal distribution}
\label{sub_dcovmvn}

In this subsection, we will make extensive use of the notation $\bN_{r \sbt}$, $\bN_{\sbt c}$,  $\bN_{* \sbt}$, $\bN_{\sbt *}$, and $\bN_{\sbt \sbt}$ from Subsection \ref{lancaster_mvn} for the multi-index matrix $\bN \in \N_0^{p \times q}$.  We now establish the following result.  

\begin{proposition}
\label{dcov_multivariatenormal}
Suppose that $(X,Y) \sim \mathcal{N}_{p+q}(0,\Sigma)$, where $\Sigma$ is given in (\ref{Sigmamatrix}).  Then the affinely invariant distance covariance, $\widetilde{\V}^2(X,Y)$, is given by 
\begin{equation}
\label{dcov_mvn}
\widetilde{\V}^2(X,Y) = \frac{1}{\gamma_p \, \gamma_q} \sum_{\bJ \neq \bzero, \bK \neq \bzero} A_{\bJ,\bK} \, B_{\bJ,\bK}  \, \frac{\Lambda_{XY}^{\bJ}}{\bJ!} \frac{\Lambda_{XY}^{\bK}}{\bK!},
\end{equation}
where the sums are taken over all non-zero $\bJ, \bK \in \N_0^{p \times q}$ such that all components of $\bJ_{* \sbts} + \bK_{* \sbts}$ and $\bJ_{\sbts *} + \bK_{\sbts *}$ are even, 
\begin{align}
\label{SbJbK}
A_{\bJ,\bK} = \frac{\Gamma\left(\tfrac12(\bJ_{\sbts \sbts} + \bK_{\sbts \sbts} - 1)\right)}{\Gamma\left(\tfrac12(\bJ_{\sbts \sbts} + \bK_{\sbts \sbts}) +\tfrac12 p\right)} \, \prod_{r=1}^p \Gamma\big(\tfrac12(\bJ_{r \sbt} + \bK_{r \sbt} + 1)\big)
\end{align}
and 
\begin{align}
\label{TbJbK}
B_{\bJ,\bK} = \frac{\Gamma\left(\tfrac12(\bJ_{\sbts \sbts} + \bK_{\sbts \sbts} - 1)\right)}{\Gamma\left(\tfrac12(\bJ_{\sbts \sbts} + \bK_{\sbts \sbts}) +\tfrac12 q\right)} \, \prod_{c=1}^q \Gamma\big(\tfrac12(\bJ_{\sbt c} + \bK_{\sbt c} + 1)\big).
\end{align}
\end{proposition}

\Pro 
In this case, the index set $\mathcal{C}$ is $\N_0^{p \times q}$, and we write the Lancaster expansion (\ref{multi.normal.lancaster}) of $(\widetilde{X},\widetilde{Y})$ in the form 
$$
\phi_{\widetilde{X},\widetilde{Y}}(x,y) - \phi_{\widetilde{X}}(x) \, \phi_{\widetilde{Y}}(y) = \phi_{\widetilde{X}}(x) \, \phi_{\widetilde{Y}}(y) \sum_{\bN \neq \bzero} \frac{\Lambda_{XY}^{\bN}}{\bN!} \, H_{\bN_{* \sbts}}(x;I_p) \, H_{\bN_{\sbts *}}(y;I_q).
$$
To calculate the Fourier transform $\P_\bN$ corresponding to $\widetilde{X}$, we apply the definition (\ref{mvhermitepoly}) of the multivariate Hermite polynomials and integration-by-parts to deduce that for $s \in \R^p$, 
\begin{align*}
\P_\bN(s) &= \int_{\R^p} \exp(\im \langle s,x\rangle) \, \phi_{\widetilde{X}}(x) \, H_{\bN_{* \sbts}}(x;I_p) \, \dd x \\
&= (-1)^{\bN_{\sbts \sbts}} \int_{\R^p} \exp(\im\langle s,x\rangle) \, \left(\frac{\partial}{\partial x}\right)^{\bN_{* \sbts}} \, \phi_{\widetilde{X}}(x) \, \dd x \\
&= \int_{\R^p} \phi_{\widetilde{X}}(x) \, \left(\frac{\partial}{\partial x}\right)^{\bN_{* \sbts}} \exp(\im\langle s,x\rangle) \, \dd x \\
&= (\im s)^{\bN_{* \sbts}} \, \int_{\R^p} \phi_{\widetilde{X}}(x) \, \exp(\im\langle s,x\rangle) \, \dd x \\
&= \im^{\bN_{\sbts\sbts}} \, s^{\bN_{* \sbts}} \, \exp(-\tfrac12\langle s,s\rangle).
\end{align*}
Similarly, 
\begin{align*}
\Q_\bN(t) = \im^{\bN_{\sbts\sbts}} \, t^{\bN_{\sbts *}} \, \exp(-\tfrac12\langle t,t\rangle),
\end{align*}
$t \in \R^q$.  Therefore, 
\begin{align*}
\int_{\R^p} \P_\bJ(s) \, \P_\bK(-s) \frac{\dd s}{\|s\|^{p+1}} = (-1)^{\bK_{\sbts \sbts}} \, \im^{\bJ_{\sbts \sbts}+\bK_{\sbts \sbts}} \, \int_{\R^p} s^{\bJ_{* \sbts}+\bK_{* \sbts}} \, \exp(-\langle s,s\rangle) \frac{\dd s}{\|s\|^{p+1}}.
\end{align*}
We now change variables to hyperspherical coordinates: $s = r \omega$, where $r > 0$ and $\omega = (\omega_1,\ldots,\omega_p) \in S^{p-1}$, the unit sphere in $\R^p$.  Then the latter integral reduces to 
$$
\int_{\R_+} r^{\bJ_{\sbts \sbts}+\bK_{\sbts \sbts}-2} \, \exp(-r^2) \, \dd r \cdot \int_{S^{p-1}} \omega^{\bJ_{* \sbts}+\bK_{* \sbts}} \dd\omega.
$$
The integral over $\R_+$ is evaluated by replacing $r$ by $r^{1/2}$, and we obtain its value as $\tfrac12 \, \Gamma\left(\tfrac12(\bJ_{\sbts \sbts}+\bK_{\sbts \sbts}-1)\right)$.  

It is easy to see that the integral over $S^{p-1}$ equals zero if any component of $\bJ_{* \sbt}+\bK_{* \sbt}$ is odd.  For the case in which each component of $\bJ_{* \sbt}+\bK_{* \sbt}$ is even, we obtain 
$$
\int_{S^{p-1}} \omega^{\bJ_{* \sbts}+\bK_{* \sbts}} \dd\omega = A(S^{p-1}) \, \E( \omega^{\bJ_{* \sbts}+\bK_{* \sbts}}),
$$
where $A(S^{p-1}) = 2 \pi^{p/2}/\Gamma(\tfrac12 p)$ is the surface area of $S^{p-1}$ and $\omega$ now is a uniformly distributed random vector on $S^{p-1}$.  It is well-known that the random vector $(\omega_1^2,\ldots,\omega_p^2) \sim D(\tfrac12,\ldots,\tfrac12)$, a Dirichlet distribution with parameters $(\tfrac12,\ldots,\tfrac12)$; so, by a classical formula for the moments of the Dirichlet distribution \cite[p.~488]{kotz2000continuous}, 
$$
\E (\omega^{\bJ_{* \sbt}+\bK_{* \sbt}}) = \frac{\Gamma(\tfrac12p)}{[\Gamma(\tfrac12)]^p} \frac{\prod_{r=1}^p \Gamma(\tfrac12(\bJ_{r \sbt} + \bK_{r \sbt} + 1))}{\Gamma(\tfrac12(\bJ_{\sbts \sbts}+\bK_{\sbts \sbts})+\tfrac12p)}.
$$
Collecting together these results, we obtain 
\begin{align*}
\int_{\R^p} \P_\bJ(s) \, \P_\bK(-s) \frac{\dd s}{\|s\|^{p+1}} = (-1)^{\bK_{\sbts \sbts}} \,(-1)^{(\bJ_{\sbts \sbts}+\bK_{\sbts \sbts})/2} \, A_{\bJ,\bK},
\end{align*}
where $A_{\bJ,\bK}$ is given in (\ref{SbJbK}).  A similar expression can be obtained for 
$$
\int_{\R^q} \Q_\bJ(t) \, \Q_\bK(-t) \frac{\dd t}{\|t\|^{q+1}}, 
$$
from which the final result (\ref{dcov_mvn}) follows. 
$\qed$

\medskip

As a consequence of Proposition \ref{dcov_multivariatenormal}, we now derive the value of the affinely invariant distance variance, $\widetilde{\V}^2(X,X)$, when $X$ has a multivariate normal distribution.  We remark that the derivation of this result given by Dueck, et al. \cite[Corollary 3.3]{dueck2014affinely} utilized the theory of zonal polynomials, whereas the proof which we now give is by a simpler method.  

\begin{corollary} {\rm (Dueck, et al. \cite[Corollary 3.3]{dueck2014affinely})}
\label{dcov_multivariatenormal_corollary}
Suppose that $(X,Y) \sim \mathcal{N}_{2p}(0,\Sigma)$, where $\Sigma$ in (\ref{Sigmamatrix}) satisfies $\Lambda_{XY} = \rho \, I_p$ in (\ref{Lambdamatrix}). Then
	\begin{equation} 
	\label{eq:cormvn}
	\widetilde{\V}^2(X,Y) = 4 \pi \frac{\gamma_{p-1}^2}{\gamma_p^2} \left[{}_2F_1(-\tfrac{1}{2},-\tfrac{1}{2};\tfrac12 p;\rho^2) - 2 \ {}_2F_1(-\tfrac{1}{2},-\tfrac{1}{2};\tfrac12 p;\tfrac14 \rho^2) +1  \right]
	\end{equation}
and 
	\begin{equation} 
	\label{eq:varmvn}
		 \widetilde{\V}^2(X,X) = 4 \pi \frac{\gamma_{p-1}^2}{\gamma_{p}^2} \, \left[
		\frac{\Gamma(\tfrac12 p) \, \Gamma(\tfrac12 p + 1)}
		{\big[\Gamma\big(\tfrac12(p+1)\big)\big]^2} 
		- 2 \ {}_2F_1 \! \left(-\tfrac12,-\tfrac12;\tfrac12 p;\tfrac14\right) + 1\right].
	\end{equation}
\end{corollary}

\Pro
Setting $p=q$ and $\Lambda_{XY} = \rho I_p$ in (\ref{dcov_mvn}), we obtain 
\begin{align}
\label{eq:mvnsepprelim}
\gamma_p^2 \, \widetilde{\V}^2(X,Y) &= \operatornamewithlimits{\sum\sum}\limits_{\substack{\bJ \neq \bzero, \bK \neq \bzero \\ \bJ_{* \sbts} + \bK_{* \sbts}, \bJ_{\sbts *} + \bK_{\sbts *} \ {\rm{even}}}} A_{\bJ,\bK} \, B_{\bJ,\bK} \, \frac{(\rho I_p)^{\bJ+\bK}}{\bJ! \, \bK!}.
\end{align}
By decomposing the set of all $(\bJ,\bK)$ into a disjoint union, 
\begin{align*}
\{(\bJ,\bK): \bJ+\bK \neq \bzero\} &= \{(\bJ,\bK): \bJ \neq \bzero, \bK \neq \bzero\} \\
& \qquad \cup \{(\bJ,\bK): \bJ = \bzero, \bK \neq \bzero\} \\
& \qquad\qquad \cup \{(\bJ,\bK): \bJ \neq \bzero, \bK = \bzero\},
\end{align*}
and noting that the summand in (\ref{eq:mvnsepprelim}) is symmetric in $(\bJ,\bK)$, we obtain 
\begin{align}
\label{eq:mvnsep}
\gamma_p^2 \, \widetilde{\V}^2(X,Y) &= \operatornamewithlimits{\sum\sum}\limits_{\substack{\bJ + \bK \neq \bzero \\ \bJ_{* \sbts} + \bK_{* \sbts}, \bJ_{\sbts *} + \bK_{\sbts *} \ {\rm{even}}}}
 A_{\bJ,\bK} \, B_{\bJ,\bK}  \, \frac{(\rho I_p)^{\bJ+\bK}}{\bJ! \, \bK!} \nonumber \\
& \qquad\qquad\qquad - 2 \operatornamewithlimits{\sum\sum}\limits_{\substack{\bJ = \bzero, \bK \neq \bzero \\ \bJ_{* \sbts} + \bK_{* \sbts}, \bJ_{\sbts *} + \bK_{\sbts *} \ {\rm{even}}}} A_{\bJ,\bK} \, B_{\bJ,\bK}  \, \frac{(\rho I_p)^{\bJ+\bK}}{\bJ! \, \bK!}.
\end{align}
Note that 
\begin{align*}
(\rho I_p)^{\bJ+\bK} = \prod_{r=1}^p \prod_{c=1}^p (\rho \delta_{rc})^{J_{rc} + K_{rc}} 
= \prod_{r=1}^p (\rho \delta_{rr})^{J_{rr} + K_{rr}} \cdot \prod_{r \neq c} (\rho \delta_{rc})^{J_{rc} + K_{rc}}.
\end{align*}
This is non-zero iff $J_{rc} = K_{rc} = 0$ for all $r \neq c$, in which case $\bJ$ and $\bK$ are diagonal matrices, and then we have 
\begin{align*}
(\rho I_p)^{\bJ+\bK} &= \rho^{J_{11}+\cdots+J_{pp}+K_{11}+\cdots+K_{pp}}, \\
\bJ_{* \sbts} + \bK_{* \sbts} &= (J_{11}+K_{11},\ldots,J_{pp}+K_{pp}) = \bJ_{\sbts *} + \bK_{\sbts *}, \\
\intertext{and }
\bJ! &= \prod_{r=1}^p J_{rr}!, \qquad  \bK! = \prod_{r=1}^p K_{rr}!.
\end{align*}
Therefore, 
$$
A_{\bJ,\bK} = B_{\bJ,\bK} = \frac{\Gamma\left(\tfrac12(J_{11}+\cdots+J_{pp}+K_{11}+\cdots+K_{pp} - 1)\right)}{\Gamma\left(\tfrac12(J_{11}+\cdots+J_{pp}+K_{11}+\cdots+K_{pp}) +\tfrac12 p\right)} \prod_{r=1}^p \Gamma\big(\tfrac12(J_{rr} + K_{rr} + 1)\big),
$$
and this yields for the first term in (\ref{eq:mvnsep})
\begin{align}
& \operatornamewithlimits{\sum\sum}\limits_{\substack{\bJ + \bK \neq \bzero \\ \bJ_{* \sbts} + \bK_{* \sbts}, \bJ_{\sbts *} + \bK_{\sbts *} \, {\rm{even}}}} A_{\bJ,\bK} \, B_{\bJ,\bK}  \, \frac{(\rho I_p)^{\bJ+\bK}}{\bJ! \, \bK!} \nonumber \\
&\qquad = \sum_{\substack{n_1,\ldots,n_p \ge 0 \\ n_1+\cdots+n_p \neq 0}} \ \sum_{J_{11}+K_{11}=2 n_1,\ldots,J_{pp}+K_{pp}=2 n_p} A_{\bJ,\bK}^2 \, \frac{\rho^{J_{11}+\cdots+J_{pp}+K_{11}+\cdots+K_{pp}}}{\bJ! \, \bK!} \nonumber \\
&\qquad = \sum_{\substack{n_1,\ldots,n_p \ge 0 \, {\rm{and \, even}} \\ n_1+\cdots+n_p \neq 0}} \rho^{2 n_1+\cdots+2 n_p} \bigg[\frac{\Gamma \left(n_1+\cdots+n_p - \tfrac12\right)}{\Gamma\left(n_1+\cdots+n_p +\tfrac12 p\right)} \, \prod_{r=1}^p \Gamma\big(n_r + \tfrac12\big)\bigg]^2 \nonumber  \\
&\qquad\qquad\qquad\qquad \times \sum_{J_{11}+K_{11}=2 n_1,\ldots,J_{pp}+K_{pp}=2 n_p} \frac{1}{\bJ! \, \bK!}. \label{eq:mvnfirst}
\end{align}
Since 
\begin{align*}
\sum_{J_{11}+K_{11}=2 n_1,\ldots,J_{pp}+K_{pp}=2 n_p} \frac{1}{\bJ! \, \bK!} = \prod_{r=1}^p \bigg[\sum_{J_{rr}+K_{rr}=2 n_r} \frac{1}{J_{rr}! \, K_{rr}!}\bigg]  
= \prod_{r=1}^p \frac{2^{2 n_r}}{(2 n_r)!},
\end{align*}
we obtain  that  (\ref{eq:mvnfirst}) equals
\begin{multline*}
\sum_{\substack{n_1,\ldots,n_p \ge 0 \\ n_1+\cdots+n_p \neq 0}} \bigg[\frac{\Gamma\left(n_1+\cdots+n_p - \tfrac12\right)}{\Gamma\left(n_1+\cdots+n_p +\tfrac12 p\right)} \, \prod_{r=1}^p \Gamma\big(n_r + \tfrac12\big)\bigg]^2 \prod_{r=1}^p \frac{2^{2n_r}\rho^{2n_r}}{(2n_r)!} \\
= \sum_{n=1}^\infty \bigg[\frac{\Gamma\left(n-\tfrac12\right)}{\Gamma\left(n+\tfrac12p\right)}\bigg]^2 (2\rho)^{2n} \sum_{n_1+\cdots+n_p=n} \prod_{r=1}^p \frac{\big[\Gamma\big(n_r + \tfrac12\big)\big]^2}{(2n_r)!}.
\end{multline*}
However, 
\begin{align*}
\sum_{n_1+\cdots+n_p=n} \prod_{r=1}^p \frac{\big[\Gamma\big(n_r + \tfrac12\big)\big]^2}{(2n_r)!} &= \sum_{n_1+\cdots+n_p=n} \prod_{r=1}^p \frac{\big[(\tfrac12)_{n_r} \Gamma(\tfrac12)\big]^2}{2^{2n_r} n_r! (\tfrac12)_{n_r}} \\
&= \frac{\pi^p}{2^{2n}} \sum_{n_1+\cdots+n_p=n} \prod_{r=1}^p \frac{(\tfrac12)_{n_r}}{n_r!} \\
&= \frac{\pi^p}{2^{2n}} \frac{(\tfrac12 p )_n}{n!},
\end{align*}
so we obtain that the first term in (\ref{eq:mvnsep}) is given by
\begin{align*}
\sum_{n=1}^\infty \bigg[\frac{\Gamma\left(n-\tfrac12\right)}{\Gamma\left(n+\tfrac12p\right)}\bigg]^2 (2\rho)^{2n} \frac{\pi^p}{2^{2n}} \frac{(\tfrac12 p)_n}{n!} 
&= \sum_{n=1}^\infty \bigg[\frac{\Gamma(-\tfrac12) (-\tfrac12)_n}{\Gamma(\tfrac12p) (\tfrac12p)_n}\bigg]^2 \rho^{2n} \pi^p \frac{(\tfrac12 p)_n}{n!} \\
&= \pi^{p} \frac{[\Gamma(-\tfrac12)]^2}{[\Gamma(\tfrac12p)]^2} \sum_{n=1}^\infty \frac{(-\tfrac12)_n (-\tfrac12)_n }{(\tfrac12p)_n} \frac{\rho^{2n}}{n!} \\
&= 4 \pi  \gamma_{p-1}^2  \left[{}_2F_1(-\tfrac12,-\tfrac12;\tfrac12p;\rho^2) - 1\right].
\end{align*}

By similar arguments, we obtain for the second term in (\ref{eq:mvnsep}): 
\begin{align*}
\operatornamewithlimits{\sum\sum}\limits_{\substack{\bJ = 0, \bK \neq \bzero \\ \bJ_{* \sbts} + \bK_{* \sbts}, \bJ_{\sbts *} + \bK_{\sbts *} \ {\rm{are \ even}}}} A_{\bJ,\bK} \, B_{\bJ,\bK}  \, \frac{(\rho I_p)^{\bJ+\bK}}{\bJ! \, \bK!} 
= 4 \pi  \gamma_{p-1}^2  \left[{}_2F_1(-\tfrac12,-\tfrac12;\tfrac12p;\tfrac14 \rho^2) - 1\right].
\end{align*}
Collecting together these results yields (\ref{eq:cormvn}).

Setting $\rho=1$ and applying Gauss' theorem for the value of ${}_2F_1(a;,b;c;1)$, we deduce that $\widetilde{\V}(X,X)$ is given by (\ref{eq:varmvn}).  
\qed

\medskip

We remark that the absolute convergence of the series in Corollary \ref{dcov_multivariatenormal_corollary} follows from the absolute convergence of Gauss' hypergeometric series.  As a consequence, the series (\ref{dcov_mvn}) converges absolutely because the matrix $\Lambda$ has norm less than $1$.

\subsection{The bivariate gamma distribution}
\label{sub_dcovbvgamma}

\begin{proposition}
\label{dcov_bivariategamma}
Suppose that the random vector $(X,Y)$ is distributed according to a Sarmanov bivariate gamma distribution, as given by (\ref{density.gamma}).  Then, 
\begin{equation}
\label{Dcov_Gamma}
\V^2(X,Y) = 2^{2(2-\alpha-\beta)}  \sum_{j, k=1}^\infty a_j \, a_k \, A_{j,k}(\alpha) \, A_{j,k}(\beta),
\end{equation}
where 
\begin{align*}
A_{j,k}(\alpha) &= 2^{-j-k} \, \left(\frac{(\alpha)_j \, (\alpha)_k}{j! \, k!} \right)^{1/2} \\
& \qquad \times \frac{\Gamma(2 \alpha+j+k-1)}{\Gamma(\alpha+j) \, \Gamma(\alpha+k)}\ {}_2F_1\left(-j-k+2,1-\alpha-j;2-2 \alpha-j-k;2 \right). 
\end{align*}
\end{proposition}

\Pro
By (\ref{density.gamma}), there holds the expansion, 
$$
\phi_{X,Y}(x,y)-\phi_X(x)\,\phi_Y(y) = \,\phi_X(x) \, \phi_Y(y) \sum_{n=1}^\infty a_n L_n^{(\alpha-1)}(x)\,L_n^{(\beta-1)}(y),
$$
$x, y > 0$.  Then, it follows from (\ref{lnx}) that for $s \in \R$, 
\begin{align*}
\P_n(s) &= \int_0^\infty \exp(\im sx) \, L_n^{(\alpha-1)}(x) \, \phi_X(x) \, \dd x \\
        &= \frac{1}{\Gamma(\alpha)} \int_0^\infty \exp\big(-(1-\im s)x\big) \, x^{\alpha-1} \, L_n^{(\alpha-1)}(x) \, \dd x.
\end{align*}
By a direct calculation using (\ref{laguerre_explicit}), we obtain 
\begin{align*}
\P_n(s) &= \left(\frac{(\alpha)_n}{n!}\right)^{1/2} \,(1-\im s)^{-\alpha} \, \left(1-(1-\im s)^{-1}\right)^n \\
&= \left(\frac{(\alpha)_n}{n!}\right)^{1/2}  \,(1-\im s)^{-(\alpha+n)} \, (-\im s)^n
\end{align*}
and, analogously,
$$
\Q_n(t) = \left(\frac{(\beta)_n}{n!}\right)^{1/2}  \, (1-\im t)^{-(\beta+n)} \, (-\im t)^n,
$$
$t \in \R$.

We now calculate the integral 
\begin{equation}
\label{g_integral}
\int_{\R} \P_j(s) \, \P_k(-s) \, \frac{\dd s}{s^2} 
\equiv \left(\frac{(\alpha)_j}{j!} \, \frac{(\alpha)_k}{k!} \right)^{1/2} \, \im^{-j+k}  \int_{\R} g(s) \dd s,
\end{equation}
where 
\begin{equation}
\label{g_function}
g(s) = s^{j+k-2} \, (1-\im s)^{-(\alpha+j)} \, (1+\im s)^{-(\alpha+k)},
\end{equation}
$s \in \R$.  To calculate the integral on the right-hand side of (\ref{g_integral}), we utilize Cauchy's beta integral \cite[p.~48]{andrews1999special} which provides that, for $a,u,v \in \C$ such that $\re(a) > 0$ and $\re(u+v) > 1$, 
\begin{equation}
\label{cauchy}
\int_{\R} (1-\im s)^{-u} \, (1+\im as)^{-v} \, \dd s = 2\pi \frac{\Gamma(u+v-1)}{\Gamma(u)\Gamma(v)} \, a^{u-1} \, (a+1)^{1-u-v}.
\end{equation}

To differentiate the left-hand side of (\ref{cauchy}) $m$ times with respect to $a$, we apply the formula, 
$$
\Big(\frac{\partial}{\partial a}\Big)^m (1+\im as)^{-v} = (- \im)^m s^m (v)_m \, (1+\im as)^{-v-m};
$$
by differentiating under the integral we obtain 
$$
\Big(\frac{\partial}{\partial a}\Big)^m \int_{\R} (1-\im s)^{-u} \, (1+\im as)^{-v} \dd s = 
(- \im)^m (v)_m \, \int_{\R} s^m \, (1-\im s)^{-u} \, (1+\im as)^{-v-m} \, \dd s.
$$

To differentiate the right-hand side of (\ref{cauchy}) $m$ times with respect to $a$, we apply Leibniz's formula,
$$
\Big(\frac{\partial}{\partial a}\Big)^m \Big[a^{u-1} \, (a+1)^{1-u-v}\Big] = \sum_{\ell=0}^m \binom{m}{\ell} \Big[\Big(\frac{\partial}{\partial a}\Big)^\ell a^{u-1}\Big] \cdot \, \Big[ \Big(\frac{\partial}{\partial a}\Big)^{m-\ell} (a+1)^{1-u-v}\Big].
$$
Noting that 
\begin{align*}
\binom{m}{\ell} &= \frac{(-1)^\ell (-m)_\ell}{\ell!}, \\
\Big(\frac{\partial}{\partial a}\Big)^{\ell} a^{u-1} &= (-1)^\ell \, a^{u-1-\ell} \, (1-u)_{\ell}, \\
\intertext{and}
\Big(\frac{\partial}{\partial a}\Big)^{m-\ell} (a+1)^{1-u-v} &= (-1)^m \, (a+1)^{1-u-v-m+\ell} \, \frac{(u+v-1)_m}{(2-u-v-m)_\ell},
\end{align*}
we obtain 
\begin{align*}
&\Big(\frac{\partial}{\partial a}\Big)^m \Big[a^{u-1} \, (a+1)^{1-u-v}\Big] \\
&= (-1)^m a^{u-1} (a+1)^{1-u-v-m} (u+v-1)_m \, \sum_{\ell=0}^m \frac{(-m)_\ell \, (1-u)_\ell}{\ell! \, (2-u-v-m)_\ell} \, a^{-\ell} \, (a+1)^{\ell} \\
&=  (-1)^m a^{u-1} (a+1)^{1-u-v-m} (u+v-1)_m \ {}_2F_1\Big(-m,1-u;2-u-v-m;\frac{a+1}{a}\Big).
\end{align*}
Comparing the derivatives of the left- and right-hand sides of (\ref{cauchy}), we obtain 
\begin{align*} 
\int_{\R} s^m \, (1-\im s)^{-u} & \, (1+\im as)^{-v-m} \, \dd s \\
&= 2\pi \, (-\im)^m \, a^{u-1} (a+1)^{1-u-v-m} \, \frac{\Gamma(u+v-1)}{\Gamma(u)\Gamma(v)} \\
& \qquad \times \frac{(u+v-1)_m}{(v)_m} \, {}_2F_1\Big(-m,1-u;2-u-v-m;\frac{a+1}{a}\Big).
\end{align*}
Substituting $a = 1$, $m = j+k-2$, $u = \alpha+j$, and $v = \alpha+k-m \equiv \alpha-j+2$, the latter equation reduces to 
\begin{align*}
\int_\R g(s) \dd s = \ & 2^{-2\alpha-j-k+2} \, \pi \, (-\im)^{j+k-2} \, \frac{\Gamma(2 \alpha+1)}{\Gamma(\alpha+j)\Gamma(\alpha-j+2)} \\ 
& \times \frac{(2 \alpha +1)_{j+k-2}}{(\alpha-j+2)_{j+k-2}} \ {}_2F_1\big(-j-k+2,1-\alpha-j;2-2 \alpha-j-k;2\big).
\end{align*}
Therefore, 
\begin{align*}
\int_{\R} \P_j(s) \, \P_k(-s) \, \frac{\dd s}{s^2} = \ & 2^{-2\alpha+1} \pi \, (-1)^{j-1} \, \left( \frac{(\alpha)_j \, (\alpha)_k}{j! \, k!} \right)^{1/2} \, \frac{\Gamma(2\alpha+j+k-1)}{\Gamma(\alpha+j)\Gamma(\alpha+k)} \\
& \times {}_2F_1\big(-j-k+2,1-\alpha-j;2-2 \alpha-j-k;2\big),
\end{align*}
and similarly for $Y$.  Substituting these expressions into Theorem \ref{Dcov_Lancaster.theorem} and simplifying the outcome, we obtain the series (\ref{Dcov_Gamma}) as a formal expression for $\mathcal{V}^2(X,Y)$.  

Finally, we verify that (\ref{Dcov_Gamma}) converges absolutely.  By (\ref{g_function}), 
\begin{align*}
\int_\R |g(s)| \, \dd s &= \int_\R |s|^{j+k-2} \, (1+s^2)^{-(2\alpha+j+k)/2} \, \dd s.
\end{align*}
Making the change-of-variables $s^2=t/(1-t)$, the latter integral is transformed to 
\begin{equation}
\label{betaintegral}
\int_0^1 t^{\frac12(j+k-3)} \, (1-t)^{\alpha-\frac12} \, \dd t = B\left(\tfrac12(j+k-1),\alpha+\tfrac12\right),
\end{equation}
where $B(\cdot,\cdot)$ is the classical beta function, and this integral converges absolutely because $j+k-1 > 0$ and $\alpha + 1/2 > 0$ for all $j, k \in \N$ and $\alpha > 0$.  Hence, to establish that (\ref{Dcov_Gamma}) converges absolutely, we have only to show that the series 
\begin{align} 
\label{double_sum}
\sum_{j=1}^\infty \sum_{k=1}^\infty a_j \, a_k \, & \left(\frac{(\alpha)_j (\beta)_j}{(j!)^2} \right)^{1/2} \, \left( \frac{(\alpha)_k \, (\beta)_k}{(k!)^2} \right)^{1/2 }\nonumber \\
& \times B\left(\tfrac12 (j+k-1),\alpha+\tfrac12 \right) \, B\left(\tfrac12(j+k-1),\beta+\tfrac12 \right)
\end{align}
converges absolutely.  

For $j+k \ge 3$, it follows from (\ref{betaintegral}) that 
$$
B\left(\tfrac12 (j+k-1),\alpha+\tfrac12 \right) \le \int_0^1 (1-t)^{\alpha -\frac12} \, \dd t = \frac{1}{\alpha+\frac12}.
$$
Therefore, (\ref{double_sum}) is bounded above by 
\begin{align*}
\beta^2 & \lambda^2 \, B(\tfrac12,\alpha+\tfrac12) \, B(\tfrac12,\beta+\tfrac12) + \frac{1}{(\alpha+\tfrac12)(\beta+\tfrac12)} \sum_{\substack{j, \, k \ge 1 \\ j+k \ge 3}} \frac{(\beta)_j}{(j!)} \frac{(\beta)_k}{(k!)} \lambda^{j+k} \\
& \le \beta^2 \lambda^2 \, B(\tfrac12,\alpha+\tfrac12) \, B(\tfrac12,\beta+\tfrac12) + \frac{1}{(\alpha+\tfrac12)(\beta+\tfrac12)} \Big(\sum_{j=0}^\infty \frac{(\beta)_j}{(j!)}  \lambda^j \Big) \Big(\sum_{k=0}^\infty \frac{ (\beta)_k}{(k!)}  \lambda^k \Big) \\
& = \beta^2 \lambda^2 \, B(\tfrac12,\alpha+\tfrac12) \, B(\tfrac12,\beta+\tfrac12) + \frac{1}{(\alpha+\tfrac12)(\beta+\tfrac12)} \ (1-\lambda)^{-2\beta},
\end{align*}
for $\lambda \in [0,1)$.  Hence (\ref{double_sum}), and then also (\ref{Dcov_Gamma}), converges absolutely for all $\alpha$, $\beta$ and all $\lambda \in [0,1)$.  
$\qed$

\bigskip

In calculating the distance variances $\V(X,X)$ and $\V(Y,Y)$, it is only the marginal distributions of $X$ and $Y$ which are relevant.  Therefore, we may assume that $X$ and $Y$ have any joint distribution for which the marginal distributions are gamma with parameters $\alpha$ and $\beta$, respectively.  

Holding $\alpha$ fixed and setting $\beta = \alpha$, the Sarmanov bivariate gamma distribution reduces to the Kibble-Moran distribution, and the characteristic function of $(X,Y)$ is 
\begin{equation}
\label{cf.bvgamma}
\psi_{X,Y}(t_1,t_2) = \big[(1 - \im t_1)(1 - \im t_2) + \lambda t_1t_2\big]^{-\alpha},
\end{equation}
see \cite[p. 436]{kotz2000continuous}.  Next, we let $\lambda \to 1-$; then $\psi_{X,Y}(t_1,t_2)$ converges to 
$$
\big[1 - \im (t_1+t_2)\big]^{-\alpha} \equiv \E \exp\big(\im (t_1+t_2)X\big),
$$
proving that if $\lambda = 1$ then $X = Y$, almost surely.  Therefore, the distance variance $\V(X,X)$ is a limiting case of $\V(X,Y)$, viz., 
\begin{align*} 
\V^2(X,X) &= \frac{1}{\gamma_1^2} \int_{\R^2} |\psi_{X}(s+t) - \psi_X(s) \psi_X(t)|^2 \frac{\dd s}{s^2} \frac{\dd t}{t^2} \\
&= \lim_{\lambda \to 1-} \frac{1}{\gamma_1^2} \int_{\R^2} |\psi_{X,Y}(s,t) - \psi_X(s) \psi_Y(t)|^2 \frac{\dd s}{s^2} \frac{\dd t}{t^2}\bigg|_{\beta = \alpha} \\
&= \lim_{\lambda \to 1-} \V^2(X,Y)\Big|_{\beta = \alpha}.
\end{align*} 
Analogously, by holding $\beta$ fixed and then setting $\alpha = \beta$, we obtain 
$$
\V^2(Y,Y) = \lim_{\lambda \to 1-} \V^2(X,Y)\Big|_{\alpha = \beta}.
$$

As a remark on the gamma distributions, note that if we replace $a_j$ in (\ref{Dcov_Gamma}) by $[(\beta)_j/(\alpha)_j]^{1/2} \,\lambda_j$ then the result above generalizes to the distribution functions introduced by Griffiths \cite{griffiths1969canonical}; see also (\ref{density.gamma.gen}) and (\ref{griff.coef}).

As noted in Subsection \ref{sub_dcovbvgamma}, if $\alpha \neq \beta$ then $\hbox{Corr}(X,Y) \not\equiv \lambda$, and then it is impractical to compare $\hbox{Corr}(X,Y)$ with $\rho$, the correlation coefficient in the bivariate normal case.  If $\alpha = \beta$ then $\hbox{Corr}(X,Y) = \lambda$, so we will consider only the case in which $\alpha = \beta$.  

\begin{figure}[!ht]
\label{dcor_graphs}
\captionsetup{width=0.9\textwidth}
\centering
\includegraphics[width=0.72\textwidth]{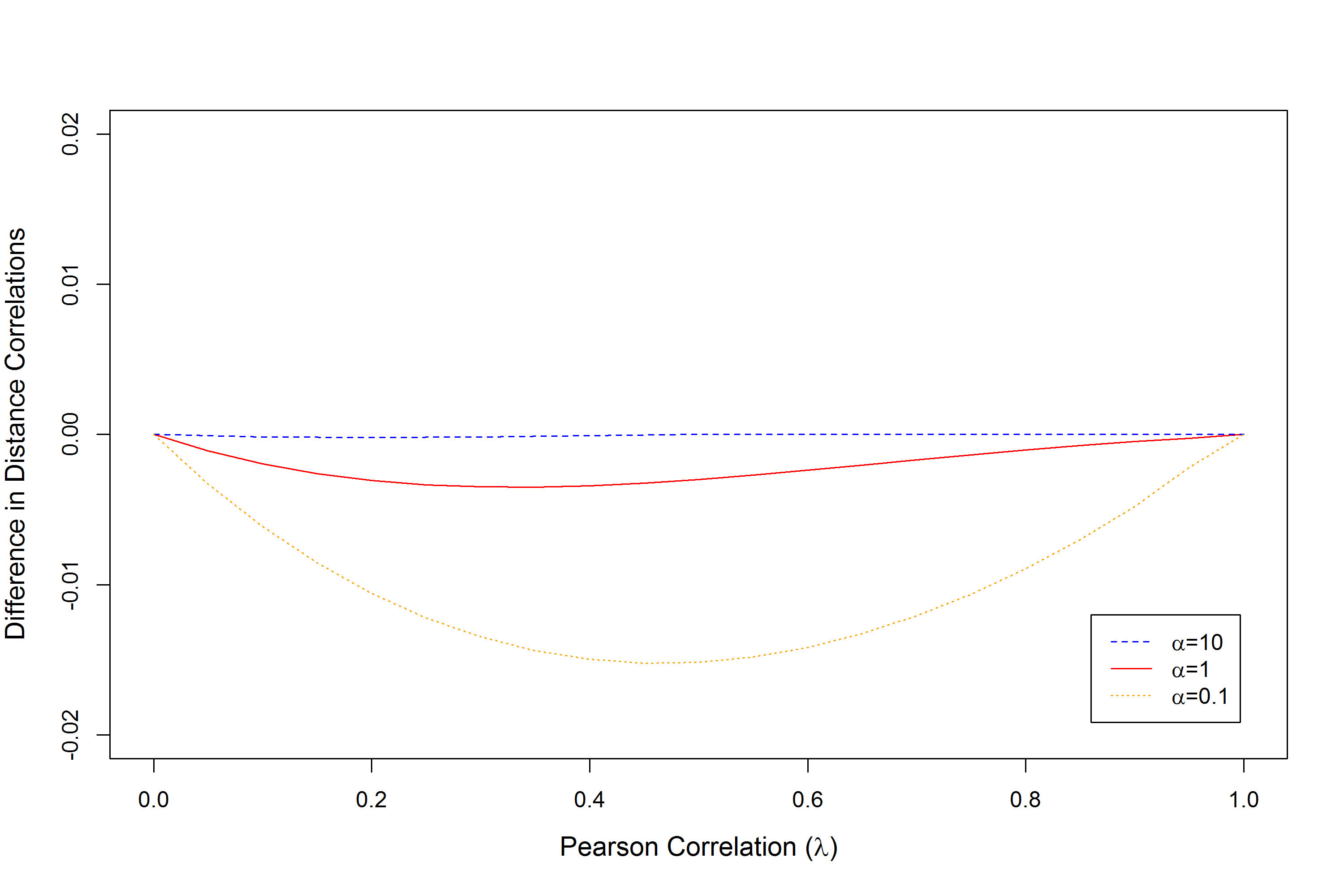}
\caption{Graphs of the difference between the distance correlation coefficient of the Kibble-Moran bivariate gamma distribution, with shape parameter $\alpha$, and the distance correlation coefficient of the bivariate normal distribution.  The graphs are given for the cases in which $\alpha = 0.1$, $1$, and $10$.}
\end{figure}

In Figure 1, we graph the difference between the distance correlation coefficient of the Kibble-Moran bivariate gamma distribution, with shape parameter $\alpha$, and the distance correlation coefficient of the bivariate normal distribution.  The graphs are given for the cases in which $\alpha = 0.1, 1, 10$.  Figure~1 suggests that the distance correlation coefficient $\mathcal{R}(X,Y)$ converges, as $\alpha \to \infty$, to the distance correlation coefficient for the bivariate normal distribution.  

This result can be proved as follows: Let $(X,Y)_\alpha$ denote a bivariate gamma random variable with probability density function (\ref{density.gamma}) where $\alpha = \beta$.  It follows from the characteristic function (\ref{cf.bvgamma}) that if $(X_1,Y_1)_\alpha,\ldots,(X_n,Y_n)_\alpha$ are independent, identically distributed random vectors with the same distribution as $(X,Y)_\alpha$ then $(X_1,Y_1)_\alpha+\cdots+(X_n,Y_n)_\alpha$ has the same distribution as $(X,Y)_{n\alpha}$.  Since $(X,Y)_\alpha$ has finite mean and covariance matrix then, by the Central Limit Theorem, $(X,Y)_{n\alpha} - \E(X,Y)_{n\alpha}$ converges, as $n \to \infty$, to a bivariate normal distribution, so $\mathcal{R}((X,Y)_{n\alpha})$ converges to the distance correlation coefficient of the bivariate normal distribution.  

Equivalently, as $\alpha \to \infty$, $\mathcal{R}((X,Y)_\alpha)$ converges to the distance correlation coefficient of the bivariate normal distribution.  Moreover, as the graphs indicates, the rate of convergence is rapid.  Indeed, for $\alpha = 0.1$, we observe from Figure 1 that the maximum absolute difference between $\mathcal{R}((X,Y)_\alpha)$ and the distance correlation of the bivariate normal distribution is less than $0.02$; and at $\alpha = 10$, the maximum absolute difference already is negligible.

\subsection{The bivariate Poisson distribution}
\label{sub_dcovbvpoisson}

\begin{proposition}
\label{dcov_bivariatepoisson}
Suppose that the random vector $(X,Y)$ is distributed according to a bivariate Poisson distribution, as given by (\ref{bivariatepoisson_pdf}).  Then 
\begin{equation}
\label{Dcov_Poisson}
\V^2(X,Y) = \sum_{j, k=1}^\infty \frac{4^{j+k-1} }{j! \, k!} \, (\lambda a)^{j+k} \, A_{jk}^2, 
\end{equation}
where 
\begin{align*}
A_{jk} = \frac{1}{(j-1)!} \sum_{\ell=0}^{\lfloor (j-k)/2\rfloor} \binom{j-k}{2\ell} (-1)^\ell (\tfrac12)_\ell \, (\tfrac12)_{j-\ell-1} \ {}_1F_1(j-\ell-\tfrac12;j;-4a)
\end{align*}
for $j \ge k$, and $A_{jk} = A_{kj}$ for $j < k$.
\end{proposition}

\Pro
By (\ref{bivariatepoisson_pdf}) and (\ref{lnx}), we have 
$$
\P_n(s) = \Q_n(s) = \E \exp(\im sX) \, C_n(X;a),
$$
$s \in \R$.  Substituting the definition (\ref{poisson_charlier}) of the Poisson-Charlier polynomials $C_n$ into the expectation and reversing the order of summation, we obtain 
\begin{align*}
\P_n(s) = \Q_n(s) &= \sum_{x=0}^\infty \exp(\im sx) C_n(x;a) \frac{\e^{-a} a^x}{x!} \\
&= \Big(\frac{a^n}{n!}\Big)^{1/2} (1 - \e^{\im s})^n \, \exp\left(-a(1-\e^{\im s})\right).
\end{align*}
Therefore, for $j, k \ge 1$, 
\begin{align}
\label{ljlkpoisson}
\int_\R & \P_j(s) \P_k(-s) \frac{\dd s}{s^2} \nonumber \\
&= \Big(\frac{a^{j+k}}{j! \, k!}\Big)^{1/2} \int_\R (1 - \e^{\im s})^j \, (1 - \e^{-\im s})^k \, \exp\left(-a(1-\e^{\im s} + 1-\e^{-\im s})\right) \frac{\dd s}{s^2} \nonumber \\
&= \Big(\frac{a^{j+k}}{j! \, k!}\Big)^{1/2} \int_\R (1 - \e^{\im s})^j \, (1 - \e^{-\im s})^k \, \exp\big(-2a(1-\cos s)\big) \frac{\dd s}{s^2}.
\end{align}
Changing variables in this integral from $s$ to $-s$ shows that the integral is symmetric in $j$ and $k$; therefore we assume, without loss of generality, that $j \ge k$.  We now write 
\begin{align*}
(1 - \e^{\im s})^j (1 - \e^{-\im s})^k &= (1 - \e^{\im s})^{j-k} (1 - \e^{\im s})^k (1 - \e^{-\im s})^k \nonumber \\
&= (1 - \e^{\im s})^{j-k} (2(1 - \cos s))^k,
\end{align*}
and apply the binomial theorem in the form, 
\begin{align*}
(1 - \e^{\im s})^{j-k} &= (1 - \cos s - \im \sin s)^{j-k} \\
&= \sum_{\ell=0}^{j-k} \binom{j-k}{\ell} (-\im \sin s)^\ell (1 - \cos s)^{j-k-\ell}.
\end{align*}
Then, it follows that the integral in (\ref{ljlkpoisson}) equals 
\begin{align}
\label{ljlkintegral}
2^k \sum_{\ell=0}^{j-k} \binom{j-k}{\ell} (-\im)^\ell \int_\R \, (\sin s)^\ell (1 - \cos s)^{j-\ell} \, \exp\big(-2a(1-\cos s)\big) \frac{\dd s}{s^2}.
\end{align}
Expanding the exponential term, 
$$
\exp\big(-2a(1-\cos s)\big) = \sum_{m=0}^\infty \frac{(-2a)^m}{m!} (1 - \cos s)^m,
$$
applying the half-angle identities, $\sin s = 2 \sin \tfrac12s \cos \tfrac12s$ and $1-\cos s = 2 (\sin \tfrac12s)^2$, and integrating term-by-term, we deduce that (\ref{ljlkintegral}) equals 
\begin{multline}
\label{ljlkintegral2}
2^k \sum_{\ell=0}^{j-k} \binom{j-k}{\ell} (-\im)^\ell \sum_{m=0}^\infty \frac{(-2a)^m}{m!} \int_\R \, (2 \sin \tfrac12s \cos \tfrac12s)^\ell (2 (\sin \tfrac12s)^2)^{j-\ell+m} \frac{\dd s}{s^2} \\
= \sum_{\ell=0}^{j-k} \binom{j-k}{\ell} (-\im)^\ell \sum_{m=0}^\infty \frac{(-a)^m}{m!} 2^{j+k+2m} \int_\R \, (\cos \tfrac12s)^\ell (\sin \tfrac12s)^{2(j+m)-\ell} \frac{\dd s}{s^2}.
\end{multline}
If $\ell$ is odd then the latter integral is an odd function of $s$, so the integral equals $0$.  Hence, (\ref{ljlkintegral2}) equals 
$$
\sum_{\ell=0}^{\lfloor (j-k)/2\rfloor} \binom{j-k}{2\ell} (-\im)^{2\ell} \sum_{m=0}^\infty \frac{(-a)^m}{m!} 2^{j+k+2m} \int_\R \, (\cos \tfrac12s)^{2\ell} (\sin \tfrac12s)^{2(j+m-\ell)} \frac{\dd s}{s^2},
$$
where $\lfloor (j-k)/2\rfloor$ denotes the greatest integer less than or equal to $(j-k)/2$.  

Next, we introduce the formula 
\begin{equation}
\label{gradshteyn_extended}
\int_\R (\cos \tfrac12s)^{2\ell} (\sin \tfrac12s)^{2k} \frac{\dd s}{s^2} = \dfrac{(2\ell-1)!! \, (2k-3)!!}{(2\ell+2k-2)!!} \, \frac{\pi}{2},
\end{equation}
$\ell=0,1,2\ldots$, $k=1,2,3,\ldots$.  This result is well-known for the case in which $\ell = 0$ (see \cite[p.~483,~3.821(10)]{gradshteyn1994table}), and the general case can be established by induction on $\ell$ with the inductive step being obtained by writing 
$$
(\cos \tfrac12s)^{2(\ell+1)} \equiv (\cos \tfrac12s)^{2\ell} (\cos \tfrac12s)^2 = (\cos \tfrac12s)^{2\ell} (1 - \sin^2 \tfrac12s).
$$
Hence, we find that (\ref{ljlkintegral2}) equals 
$$
\frac{\pi}{2} \sum_{\ell=0}^{\lfloor (j-k)/2\rfloor} \binom{j-k}{2\ell} (-1)^\ell (2\ell-1)!! \, \sum_{m=0}^\infty \frac{(-a)^m}{m!} 2^{j+k+2m} \dfrac{(2(j+m-\ell)-3)!!}{(2(j+m)-2)!!}.
$$
Writing each double factorial in terms of rising factorials, and simplifying the resulting expressions, we find that this sum equals 
\begin{multline*}
\frac{\pi}{2} \sum_{\ell=0}^{\lfloor (j-k)/2\rfloor} \binom{j-k}{2\ell} (-1)^\ell (2\ell-1)!! \, 2^{j+k-\ell} \frac{(\frac12)_{j-\ell-1}}{(j-1)!} \sum_{m=0}^\infty \frac{(-a)^m}{m!} 2^{2m} \frac{(j-\ell-\frac12)_m}{(j)_m} \\
= \pi \frac{2^{j+k-1}}{(j-1)!} \sum_{\ell=0}^{\lfloor (j-k)/2\rfloor} \binom{j-k}{2\ell} (-1)^\ell (\tfrac12)_\ell \, (\tfrac12)_{j-\ell-1} \ {}_1F_1(j-\ell-\tfrac12;j;-4a).
\end{multline*}
Substituting this result into (\ref{ljlkintegral2}), we obtain 
\begin{align*}
\int_\R & \P_j(s) \P_k(-s) \frac{\dd s}{s^2} \nonumber \\
&= \pi \Big(\frac{a^{j+k}}{j! \, k!}\Big)^{1/2} \frac{2^{j+k-1}}{(j-1)!} \sum_{\ell=0}^{\lfloor (j-k)/2\rfloor} \binom{j-k}{2\ell} (-1)^\ell (\tfrac12)_\ell \, (\tfrac12)_{j-\ell-1} \ {}_1F_1(j-\ell-\tfrac12;j;-4a).
\end{align*}
Substituting this result into Theorem \ref{Dcov_Lancaster.theorem} and simplifying the outcome, we obtain the series (\ref{Dcov_Poisson}) as a formal expression for $\V^2(X,Y)$.
 
Finally, we establish the absolute convergence of the series (\ref{Dcov_Poisson}).  On applying to (\ref{ljlkpoisson}) the identity 
$$
|1 - \e^{\im s}| = |1 - \e^{-\im s}| = \big(2(1 - \cos s)\big)^{1/2} = 2\big(\sin^2 \tfrac12 s\big)^{1/2}
$$
and the inequality 
$$
\exp\big(-2(1 - \cos s)\big) \le 1,
$$
$s \in \R$, we obtain 
\begin{align*}
\left|\int_\R \P_j(s) \P_k(-s) \frac{\dd s}{s^2}\right| &\le \Big(\frac{a^{j+k}}{j! \, k!}\Big)^{1/2} \int_\R |1 - \e^{\im s}|^j |1 - \e^{\im s}|^k \exp\big(-2a(1-\cos s)\big) \frac{\dd s}{s^2} \\
&\le \Big(\frac{(4a)^{j+k}}{j! \, k!}\Big)^{1/2} \int_\R \big(\sin^2 \tfrac12s\big)^{(j+k)/2} \frac{\dd s}{s^2}.
\end{align*}
By the Cauchy-Schwarz inequality, 
\begin{align*}
\int_\R \big(\sin^2 \tfrac12s\big)^{(j+k)/2} \frac{\dd s}{s^2} &\equiv \int_\R \big(\sin^2 \tfrac12s\big)^{j/2} \big(\sin^2 \tfrac12s\big)^{k/2} \frac{\dd s}{s^2} \\
&\le \left(\int_\R \big(\sin^2 \tfrac12s\big)^j \frac{\dd s}{s^2}\right)^{1/2} \ \left(\int_\R \big(\sin^2 \tfrac12s\big)^k \frac{\dd s}{s^2}\right)^{1/2}.
\end{align*}
Since $(2k-3)!!/(2k-2)!! \le 1$ for all $k \in \N$ then it follows from (\ref{gradshteyn_extended}) with $\ell = 0$ that 
$$
\int_\R \big(\sin^2 \tfrac12s\big)^j \frac{\dd s}{s^2} \le \pi;
$$
therefore, 
$$
\left|\int_\R \P_j(s) \P_k(-s) \frac{\dd s}{s^2}\right| \le \Big(\frac{(4a)^{j+k}}{j! \, k!}\Big)^{1/2} \, \pi,
$$
and the same holds for the functions $\Q_j$.  Substituting these bounds into the general series expansion (\ref{Dcov_Lancaster}), we obtain the upper bound 
\begin{align*}
\V^2(X,Y) \le \sum_{j=1}^\infty \sum_{k=1}^\infty \frac{(4\lambda a)^{j+k}}{j! \, k!} 
= \big(\exp(4\lambda a)-1\big)^2 < \infty,
\end{align*}
for all $\lambda \in [0,1]$ and $a > 0$.  Therefore, the series (\ref{Dcov_Poisson}) converges absolutely.  
$\qed$

\bigskip

To calculate the distance variance, the argument given in the bivariate gamma case remains valid here.  By Koudou \cite[p.~103]{koudou1998lancaster}, the characteristic function of $(X,Y)$ is 
$$
\psi_{X,Y}(s,t) = \exp\big[a(1-\lambda)(\e^{\im s} - 1) + a(1-\lambda)(\e^{\im t} - 1) + a\lambda(\e^{\im (s+t)} - 1)\big].
$$
Therefore, 
$$
\lim_{\lambda \to 1-} \psi_{X,Y}(s,t) = \exp\left[a(\e^{\im (s+t)} - 1\right] \equiv \psi_X(s+t),
$$
so we obtain 
$$
\V^2(X,X) = \V^2(Y,Y) = \lim_{\lambda \to 1-} \V^2(X,Y).
$$

\subsection{The bivariate negative binomial distribution}

\begin{proposition}
\label{dcov_bivariate_nb}
Suppose that the random vector $(X,Y)$ is distributed according to a bivariate negative binomial distribution, as given by (\ref{bivariatenb_pdf}).  Then,
\begin{equation}
\label{Dcov_nb}
\V^2(X,Y) = (1-c)^{4\beta} \sum_{j,k=1}^\infty \frac{(\beta)_j \, (\beta)_k}{j! \, k!} (1+c^2)^{-2\beta-2j} 2^{2k} (\lambda c)^{j+k} A_{jk}^2,
\end{equation}
where 
\begin{align*}
A_{jk} &= \sum_{\ell_1,\ell_2=0}^{j-k} \binom{j-k}{\ell_1} \binom{j-k}{\ell_2} (-c)^{\ell_1} (-1)^{\ell_2} (|\ell_1-\ell_2|)! \sum_{\ell=0}^\infty \frac{(\beta+j)-\ell}{\ell!} \bigg(\frac{2c}{1+c^2}\bigg)^\ell \\
& \quad \times \sum_{m=0}^{|\ell_1-\ell_2|} \frac{(-2)^m (m)_{|\ell_1-\ell_2|}}{(|\ell_1-\ell_2|-m)! \, (2m)!} \frac{2^{k+m-1} \, (\tfrac12)_{k+m-1}}{(k+m-1)!} \ {}_2F_1(-\ell,k+m-\tfrac12;k+m;2),
\end{align*}
for $j \ge k$, and $A_{jk} = A_{kj}$ for $j < k$.  
\end{proposition}

\Pro
By (\ref{bivariatenb_pdf}), 
$$
\phi_{X,Y}(x,y)-\phi_X(x)\,\phi_Y(y) = \,\phi_X(x) \, \phi_Y(y) \sum_{n=1}^\infty \lambda^n M_n^{\beta,c}(x) \, M_n^{\beta,c)}(y),
$$
$x, y \in \N_0$.  Then, it follows from (\ref{lnx}) that for $s \in \R$, 
\begin{align*}
\P_n(s) = \Q_n(s) &= \E \exp(\im sX) \, M_n^{\beta,c}(X) \\
&= \sum_{x=0}^\infty \exp(\im sx) \, M_n^{\beta,c}(x) \, (1-c)^\beta \, \frac{c^x \, (\beta)_x}{x!}.
\end{align*}
Substituting the definition of the Meixner polynomial as given in (\ref{meixner}), we obtain 
\begin{align*}
\P_n(s) &= (1-c)^\beta \, \Big(\frac{c^n \,(\beta)_n}{n!}\Big)^{1/2} \sum_{x=0}^\infty \exp(\im sx) \, \frac{c^x \, (\beta)_x}{x!} \, \sum_{k=0}^n \frac{(-n)_k \, (-x)_k}{(\beta)_k \, k!} \left(1-c^{-1}\right)^k.
\end{align*}
Letting $u \equiv c\e^{\im s}$ and interchanging the order of summation, we obtain 
\begin{equation}
\label{Pnnegativebinomial}
\P_n(s) = (1-c)^\beta \, \Big(\frac{c^n \,(\beta)_n}{n!}\Big)^{1/2} \sum_{k=0}^n \frac{(-n)_k}{(\beta)_k \, k!} (1-c^{-1})^k \sum_{x=0}^\infty \frac{u^x \, (-x)_k (\beta)_x}{x!}.
\end{equation}
Noting that 
$$
(-x)_k u^x = (-1)^k u^k \left(\frac{\dd}{\dd u}\right)^k u^x,
$$
we obtain 
\begin{align*}
\sum_{x=0}^\infty \frac{u^x \, (-x)_k (\beta)_x}{x!} &= (-1)^k u^k \left(\frac{\dd}{\dd u}\right)^k \sum_{x=0}^\infty \frac{u^x \, (\beta)_x}{x!} \\
&= (-1)^k u^k \left(\frac{\dd}{\dd u}\right)^k (1-u)^{-\beta} \\
&= (-1)^k (\beta)_k u^k (1-u)^{-\beta-k}.
\end{align*}
Substituting this result into the inner summation in (\ref{Pnnegativebinomial}) and simplifying the resulting expression, we obtain 
\begin{align*}
\P_n(s) 
&= (1-c)^\beta \, \Big(\frac{c^n \,(\beta)_n}{n!}\Big)^{1/2} (1-c\e^{\im s})^{-\beta} \sum_{k=0}^n \frac{(-n)_k}{k!} \left(\frac{(1-c) \e^{\im s}}{1-c\e^{\im s}}\right)^k \\
&= (1-c)^\beta \, \Big(\frac{c^n \,(\beta)_n}{n!}\Big)^{1/2} (1-c\e^{\im s})^{-\beta} \left(1-\frac{(1-c) \e^{\im s}}{1-c\e^{\im s}}\right)^n \\
&= (1-c)^\beta \, \Big(\frac{c^n \,(\beta)_n}{n!}\Big)^{1/2} (1-c\e^{\im s})^{-\beta-n} (1-\e^{\im s})^n.
\end{align*}
Therefore, for $j, k \ge 1$, 
\begin{align}
\label{ljlknegativebinomial}
\int_\R \P_j(s) \P_k(-s) \frac{\dd s}{s^2} &= (1-c)^{2\beta} \, \Big(\frac{c^{j+k} \,(\beta)_j (\beta)_k}{j! \,k!}\Big)^{1/2} \nonumber \\
& \quad\ \times \int_\R (1 - c\e^{\im s})^{-\beta-j} (1 - c\e^{-\im s})^{-\beta-k} \, (1 - \e^{\im s})^j \, (1 - \e^{-\im s})^k \, \frac{\dd s}{s^2}. \qquad
\end{align}
Changing the variable of integration from $s$ to $-s$ shows that (\ref{ljlknegativebinomial}) is symmetric in $j$ and $k$; so we assume, without loss of generality, that $j \ge k$. 

Next, we write the integrand in (\ref{ljlknegativebinomial}) in the form 
\begin{align*}
(1 & - c\e^{\im s})^{-\beta-j} (1 - c\e^{-\im s})^{-\beta-k} \, (1 - \e^{\im s})^j \, (1 - \e^{-\im s})^k \\
&= [(1 - c\e^{\im s})(1 - c\e^{-\im s})]^{-\beta-j} [(1 - \e^{\im s})(1 - \e^{-\im s})]^k (1 - c\e^{-\im s})^{j-k} (1 - \e^{\im s})^{j-k} \\
&= (1 + c^2 - 2c \cos s)^{-\beta-j} [2(1-\cos s)]^k (1 - c\e^{-\im s})^{j-k} (1 - \e^{\im s})^{j-k} \\
&= 2^k (1 + c^2)^{-\beta-j} \Big(1-\frac{2c}{1+c^2} \cos s\Big)^{-\beta-j} (1-\cos s)^k (1 - c\e^{-\im s})^{j-k} (1 - \e^{\im s})^{j-k}.
\end{align*}
By the binomial theorem, 
$$
(1 - c\e^{-\im s})^{j-k} (1 - \e^{\im s})^{j-k} = \sum_{\ell_1,\ell_2=0}^{j-k} \binom{j-k}{\ell_1} \binom{j-k}{\ell_2} (-c)^{\ell_1} (-1)^{\ell_2} \e^{\im (\ell_2-\ell_1) s},
$$
so it follows that the integral in (\ref{ljlknegativebinomial}) is a linear combination of integrals of the form 
$$
\int_\R \Big(1-\frac{2c}{1+c^2} \cos s\Big)^{-\beta-j} \, (1-\cos s)^k \, \e^{-\im (\ell_1-\ell_2) s} \frac{\dd s}{s^2}
$$
which, after symmetrizing with $s$ replaced by $-s$, equals 
$$
\frac12 \int_\R \Big(1-\frac{2c}{1+c^2} \cos s\Big)^{-\beta-j} \, (1-\cos s)^k \, [\e^{-\im (\ell_1-\ell_2) s} + \e^{\im (\ell_1-\ell_2) s}] \frac{\dd s}{s^2}
$$
$$
= \int_\R \Big(1-\frac{2c}{1+c^2} \cos s\Big)^{-\beta-j} \, (1-\cos s)^k \, \cos\big((\ell_1-\ell_2)s\big) \frac{\dd s}{s^2}.
$$

To calculate this integral, we write 
$$
\Big(1 - \frac{2c}{1 + c^2}\cos s\Big)^{-\beta-j} = \sum_{\ell=0}^\infty \frac{(\beta+j)_\ell}{\ell!} \Big(\frac{2c}{1 + c^2} \, \cos s\Big)^\ell,
$$
and apply the Chebyshev polynomials $T_n$, given by 
$$
\cos\big((\ell_1-\ell_2)s\big) = T_{\ell_1-\ell_2}(\cos s) 
= (|\ell_1-\ell_2|)! \sum_{m=0}^{|\ell_1-\ell_2|} (-2)^m \frac{(m)_{|\ell_1-\ell_2|}}{(|\ell_1-\ell_2|-m)! \, (2m)!} \, (1-\cos s)^m;
$$
see \cite[p. 1056, 8.942(1)]{gradshteyn1994table}.  Then, we see that we need to calculate 
$$
\int_\R (\cos s)^\ell (1-\cos s)^{k+m} \frac{\dd s}{s^2}.
$$
Using the standard half-angle transformations for the cosine function and applying (\ref{gradshteyn_extended}), we obtain 
\begin{align}
\label{eq:cosinepowerintegral}
\int_\R (\cos s)^\ell \, (1 - \cos s)^{k+m} \frac{\dd s}{s^2} &= \int_\R  (1-2\sin^2 \tfrac12 s)^\ell \, (2\sin^2 \tfrac12 s)^{k+m} \frac{\dd s}{s^2} \nonumber \\
&= 2^{k+m} \sum_{r=0}^\ell \binom{\ell}{r} (-2)^r \int_\R (\sin^2 \tfrac12 s)^{r+k+m} \frac{\dd s}{s^2} \nonumber \\
&= 2^{k+m-1} \pi \sum_{r=0}^\ell \binom{\ell}{r} (-2)^r \frac{(2r+2k+2m-3)!!}{(2r+2k+2m-2)!!}. 
\end{align}
Expressing these double factorials in terms of rising factorials, we find that (\ref{eq:cosinepowerintegral}) equals 
$$
\pi \frac{2^{k+m-1} (\tfrac12)_{k+m-1}}{(k+m-1)!} \ {}_2F_1(-\ell,k+m-\tfrac12;k+m;2).
$$
Collecting together all terms, we obtain (\ref{Dcov_nb}).

Finally, a proof of the absolute convergence of (\ref{Dcov_nb}) can be obtained using arguments similar to those used to establish convergence in the previous subsections.  
$\qed$

\section{Summary and conclusions}

We have derived a general series expansion for the distance covariance and distance correlation coefficients for the class of Lancaster distributions.  This result resolves the fundamental obstacle arising in calculating the singular integrals used to define distance correlation.  We have established the utility of the result by applying it to derive the distance correlation for the bivariate normal distribution and its generalizations of Lancaster type, the multivariate normal distributions, and the bivariate gamma, Poisson, and negative binomial distributions which are of Lancaster type.  

In computing any of the series obtained in this paper, we can derive upper bounds on the maximum discrepancy arising from the use of a finite number of terms of the series by applying well-known methods of Kotz, et al. \cite{kotz1967series} together with classical bounds for the various hypergeometric series appearing in the expansions.

\section*{Acknowledgments}

The research of Dueck and Edelmann was supported by the {\it Deutsche Forschungsgemeinschaft} (German Research Foundation) within the programme ``Spatio/Temporal Graphical Models and Applications in Image Analysis,'' grant GRK 1653.  The research of Richards was partially supported by the U.S. National Science Foundation, grant DMS-1309808; by a sabbatical leave-of-absence at Heidelberg University; and by a Romberg Guest Professorship at the Heidelberg University Graduate School for Mathematical and Computational Methods in the Sciences, funded by the German Universities Excellence Initiative grant GSC 220/2.


\begin{thebibliography}{99}

\bibitem{andrews1999special}
Andrews, G.~E., Askey, R., and Roy, R.~(1999).  {\sl Special Functions}.  
Cambridge University Press, New York.

\bibitem{bar1994diagnonal}
Bar-Lev, S.~K., Bshouty, D., Letac, G., Lu, I.,~and Richards, D.~St.~P.~(1993).
The diagonal multivariate natural exponential families and their
classification. {\em J. Theoret. Probab.}, {\bf 7}, 883--929.

\bibitem{diaconis2008gibbs}
Diaconis, P., Khare, K., and Saloff-Coste, L.~(2008). Gibbs sampling, exponential 
families and orthogonal polynomials. {\em Statist. Sci.}, {\bf 23}, 151--178.

\bibitem{dueck2012affinely}
Dueck, J., Edelmann, D., Gneiting, T., and Richards, D.~St.~P.~(2012). The affinely
invariant distance correlation. Preprint, \url{http://arxiv.org/abs/1210.2482v1}.  

\bibitem{dueck2014affinely}
Dueck, J., Edelmann, D., Gneiting, T., and Richards, D.~St.~P.~(2014). The affinely
invariant distance correlation. {\em Bernoulli}, {\bf 20}, 2305--2330.

\bibitem{dueck2015generalization}
Dueck, J., Edelmann, D., and Richards, D.~(2015). A generalization of an 
integral arising in the theory of distance correlation. 
{\em Statist. \& Probab. Lett.}, {\bf 97}, 116--119. 

\bibitem{gradshteyn1994table}
Gradshteyn, I. S., and Ryzhik, I. M. (1994).  {\sl Tables of Integrals, Series, 
and Products}, fifth edition (A. Jeffrey, Editor).  Academic Press, New York.

\bibitem{griffiths1969canonical}
Griffiths, R. C. (1969). The canonical correlation coefficients of bivariate 
gamma distributions. {\em Ann. Math. Statist.}, {\bf 40}, 1401-1408.

\bibitem{huo2015fast}
Huo, X., and Sz\'ekely, G. J. (2016).  Fast computing for distance covariance.  
{\em Technometrics}, {\bf 58}, 435--447.

\bibitem{kong2012using}
Kong,~J., Klein,~B.~E.~K., Klein,~R., and Wahba,~G.~(2012).  Using distance 
correlation and SS-ANOVA to access associations of familial relationships, 
lifestyle factors, diseases, and mortality.  {\it Proc. Natl. Acad. Sci. USA}, 
{\bf 109}, 20352--20357.

\bibitem{kotz2000continuous}
Kotz, S., Balakrishnan, N., and Johnson, N.~L.~(2000). {\em Continuous
Multivariate Distributions, Volume 1, Models and Applications}, second edition. 
Wiley, New York.

\bibitem{kotz1967series}
Kotz, S., Johnson, N. L., and Boyd, D. W. (1967).  Series representations of 
distributions of quadratic forms in normal variables. I. Central case.  
{\it Ann.Math. Statist.}, {\bf 38}, 823--837.

\bibitem{koudou1996probabilites}
Koudou, A.~E.~(1996).  Probabilit\'es de Lancaster. {\em Expo. Math.}, 
{\bf 14}, 247--275.

\bibitem{koudou1998lancaster}
Koudou, A.~E.~(1998).  Lancaster bivariate probability distributions 
with Poisson, negative binomial and gamma margins. {\em Test}, {\bf 7},
95--110.

\bibitem{lancaster1958structure}
Lancaster, H.~O.~(1958). The structure of bivariate distributions. {\em 
Ann. Math. Statist.}, {\bf 29}, 719--736.

\bibitem{lancaster1969}
Lancaster, H.~O.~(1969). {\em The Chi-Squared Distribution}. Wiley, New York.

\bibitem{lyons2013distance}
Lyons, R. (2013). Distance covariance in metric spaces.  
{\em Ann. Probab.}, {\bf 41}, 3284--3305.

\bibitem{martinez2014distance}
Mart{\'i}nez-G\'omez, E., Richards, M. T., and Richards, D. St. P.~(2014). 
Distance correlation methods for discovering associations in large 
astrophysical databases. {\em Astrophys. J.}, {\bf 781}, 39 (11 pp.).

\bibitem{pommeret2004characterization}
Pommeret, D.~(2004). A characterization of Lancaster probabilities with margins
in a multivariate additive class. {\em Sankhy\=a}, {\bf 66}, 1--19.

\bibitem{richards2014interpreting}
Richards, M. T., Richards, D. St. P., and Mart{\'i}nez-G\'omez, E. (2014). 
Interpreting the distance correlation results for the COMBO-17 survey. 
{\em Astrophys. J. Lett.}, {\bf 784}, L34 (5 pp.).

\bibitem{rizzo2010disco}
Rizzo, M. L., and Sz\'ekely, G. J. (2010).  DISCO analysis: A nonparametric 
extension of analysis of variance. {\em Ann. Appl. Statist.}, {\bf 4}, 
1034--1055.

\bibitem{sarmanov1970approximate}
Sarmanov, I.~O.~(1970).  The approximate computation of the coefficient 
of correlation between functions of dependent random variables. 
{\em Math. Notes Acad. Sciences USSR}, {\bf 7}, 373--377.

\bibitem{sarmanov1966generalized}
Sarmanov, O.~V.~(1966).  Generalized normal correlation and two-dimensional
Fr\'echet classes. {\em Soviet Math. Dokl.}, {\bf 7}, 596--599.

\bibitem{sarmanov1970gamma}
Sarmanov, O.~V.~(1970).  A gamma-correlation process and its properties. 
{\em Dokl. Akad. Nauk SSSR}, {\bf 191}, 30--32.

\bibitem{sarmanov1967probabilistic}
Sarmanov, O.~V., and Bratoeva, Z.~N.~(1967).  Probabilistic properties of 
bilinear expansions in Hermite polynomials. {\em Theor. Probab. Appl.}, 
{\bf 12}, 470--481.

\bibitem{sejdinovic2013equivalence}
Sejdinovic, D., Sriperumbudur, B., Gretton, A., and Fukumizu, K.~(2013). 
Equivalence of distance-based and RKHS-based statistics in hypothesis testing. 
{\em Ann. Statist.}, {\bf 41}, 2263--2291.


\bibitem{szekely2009brownian}
Sz\'ekely, G. J.,~and Rizzo, M.~(2009).  Brownian distance covariance.
{\em Ann. Appl. Statist.}, {\bf 3}, 1236--1265.

\bibitem{szekely2012uniqueness}
Sz\'ekely, G. J.,~and Rizzo, M.~(2012).  On the uniqueness of distance 
correlation.  {\em Statist. \& Probab. Lett.}, {\bf 82}, 2278--2282.

\bibitem{szekely2013distance}
Sz\'ekely, G. J.,~and Rizzo, M.~(2013).  The distance correlation $t$-test 
of independence in high dimension.  {\it J. Multivariate Anal.}, {\bf 117}, 
193--213.

\bibitem{szekely2014partial}
Sz\'ekely, G. J.,~and Rizzo, M.~(2014).  Partial distance correlation with 
methods for dissimilarities.  {\it Ann. Statist.}, {\bf 42}, 2382--2412.

\bibitem{szekely2007measuring}
Sz\'ekely, G.~J., Rizzo, M.~L.,~and Bakirov, N.~K.~(2007).  Measuring
and testing dependence by correlation of distances.  {\em Ann. Statist.}, 
{\bf 35}, 2769--2794.

\bibitem{withers2010expansions}
Withers, C. S.,~and Nadarajah, S.~(2010). Expansions for the 
multivariate normal. {\em J. Multivariate Anal.}, {\bf 101}, 1311--1316.

\bibitem{zhou2012measuring}
Zhou, Z.~(2012).  Measuring nonlinear dependence in time-series, a distance 
correlation approach.  {\em J. Time Series Anal.}, {\bf 33}, 438--457.

\end{thebibliography}
\end{document}